\documentclass[twoside]{irmaems}

\usepackage{amssymb} 
\usepackage{amsmath} 
\usepackage{latexsym}
\usepackage[all]{xy}
\renewcommand\theequation{\arabic{section}.\arabic{equation}}

\theoremstyle{definition} 
\newtheorem{df}{Definition}[section]
\newtheorem{remark}[df]{Remark}

\theoremstyle{plain}      
\newtheorem{prop}{Proposition}[section]
\newtheorem{lemma}[prop]{Lemma}
\newtheorem{thm}[prop]{Theorem}
\newtheorem{cor}[prop]{Corollary}

\newcommand{\II}{I\hspace{-0.1cm}I}
\newcommand{\III}{I\hspace{-0.1cm}I\hspace{-0.1cm}I}
\newcommand{\dr}{\partial}
\newcommand{\cT}{{\mathcal T}}
\newcommand{\cM}{{\mathcal M}}
\newcommand{\cTb}{{\overline{\mathcal T}}}
\newcommand{\Grb}{{\overline{Gr}}}
\newcommand{\db}{{\overline{\partial}}}
\newcommand{\cS}{{\mathcal S}}

\newcommand{\tr}{\mbox{tr}}
\newcommand{\be}{\begin{eqnarray}}
\newcommand{\ee}{\end{eqnarray}}
\newcommand{\C}{{\mathbb C}}
\newcommand{\R}{{\mathbb R}}
\newcommand{\CP}{\mathcal{CP}}
\newcommand{\ML}{\mathcal{ML}}

\newcommand{\GF}{\mathcal{GF}}
\newcommand{\QF}{\mathcal{QF}}

\newcommand{\St}{\tilde{S}}
\newcommand{\dev}{\mathrm{dev}}
\newcommand{\Eps}{\mathrm{Eps}}

\begin{document}

\markboth{Kirill Krasnov and Jean-Marc Schlenker}{Title}
\title{The Weil-Petersson metric and the renormalized volume of hyperbolic 3-manifolds}
\author{Kirill Krasnov\thanks{Supported by an EPSRC Advanced fellowship} and 
Jean-Marc Schlenker\thanks{Partially supported by the ANR programs Repsurf (ANR-06-BLAN-0311) and 
ETTT (ANR-09-BLAN-0116-01).}
}
\address{K. K.: School of Mathematical Sciences, University of
Nottingham, Nottingham, NG7 2RD, UK
\\[4pt]
J.-M. S.: Institut de Math{\'e}matiques de Toulouse, UMR CNRS 5219 \\
Universit{\'e} Toulouse III,
31062 Toulouse cedex 9, France.} 

\maketitle


\tableofcontents

\section{Introduction}

We survey the renormalized volume of hyperbolic  3-manifolds, as a tool 
for Teichm\"uller theory, using simple differential geometry arguments to
recover results sometimes first achieved by other means. One such application is 
McMullen's quasifuchsian (or more generally Kleinian) reciprocity, for which
different arguments are proposed. Another is the fact that the renormalized
volume of quasifuchsian (or more generally geometrically finite) hyperbolic
3-manifolds provides a K\"ahler potential for the Weil-Petersson metric on
Teichm\"uller space. Yet another is the fact that the grafting map is
symplectic, which is proved using a variant of the renormalized volume 
defined for hyperbolic ends.

\subsection{Liouville theory\index{Liouville theory}\index{theory!Liouville}}

One of the early approaches to the problem of uniformization of Riemann 
surfaces was based on the so-called Liouville equation
\index{Liouville equation}\index{equation!Liouville}.
Consider a closed Riemann surface $S$ of genus $g$ and fix an arbitrary reference
metric $h_0$ in the conformal class of $S$. Then, consider a conformally
equivalent metric $h=e^{2\phi} h_0$. The condition that $h$ has constant
curvature minus one reads:
\be\label{liouv}
\Delta_0 \phi - K_0 = e^{2\phi},
\ee
where $\Delta_0$ and $K_0$ are the Laplacian and Gauss curvature of $h_0$
respectively. One then tries to solve this equation for $\phi$ and thus
find a hyperbolic metric on $S$. 

Historically, this approach to the uniformization turned out to be too difficult, and
was abandoned in favor of the one based on Fuchsian groups. More recently,
the set of ideas related to the Liouville equation came to the spotlight
due to the central role it plays in Polyakov's approach to string theory.
In the so-called non-critical string theory a very important role is
played by the {\it Liouville functional}\index{Liouville functional}\index{functional!Liouville}, 
which (in one of its versions) 
can be written as:
\be\label{liouv-1}
S[h_0,\phi]=\frac{1}{8\pi}\int d{\rm vol}_0 
\left( |\nabla \phi|^2 + e^{2\phi} - 2\phi K_0\right).
\ee


When varied with respect to $\phi$ this functional gives the Liouville 
equation (\ref{liouv}). As is implied by the uniformization theorem,
there is indeed a unique solution to (\ref{liouv}) on any given Riemann
surface $S$. Let us denote this solution by $\phi_{{\rm hyp}}$, where ${\rm hyp}$ stands for
hyperbolic. One can evaluate the functional (\ref{liouv-1}) on $\phi_{{\rm hyp}}$ and
obtain a functional $S[h_0]=S[h_0,\phi_{{\rm hyp}}]$. When $h_0$ is taken to be the 
hyperbolic metric, $\phi_{{\rm hyp}}=0$ and the value of the above functional is
just $1/4\pi$ times the area of $S$ evaluated using the hyperbolic metric, i.e.
minus half the Euler characteristics of $X$. The values of $S[h_0]$ at other points
are not easy to find and one gets a rather non-trivial functional on the space
of metrics $h_0$ on $S$, which, in a certain sense, measures how far the metric
$h_0$ is from the hyperbolic one. This functional is not so interesting
by itself, but serves as a prototype for the construction of a whole class of 
Liouville-type functionals that will play the central role in this article. The
reason why $S[h_0]$ is not so interesting is that one would rather have
a functional on the space of conformal equivalence classes of metrics on $S$,
i.e. the {\it moduli space} of Riemann surfaces, or at
least on the Teichm\"uller space $\cT_g$, hoping that such a functional could
be used to characterize the geometry of the moduli space in an interesting way.
Such functionals are not easy to construct. Indeed, one could have obtained a functional
on the moduli from $S[h_0]$ by taking $h_0$ to be some canonical metric in the given conformal
class. However, as we have already seen, taking $h_0$ to be the canonical metric of constant
negative curvature gives a functional on the moduli space whose value at all
points is a constant --- the Euler characteristic. 

The reader can find another, discrete approach to Liouville theory in the survey
by Kashaev in this volume \cite{kashaev-survey}.

\subsection{Liouville theory and projective structures\index{projective}}

It turned out to be possible to use Liouville theory for a construction of interesting functions 
on the moduli space by first considering a certain function on the space of complex projective structures 
on Riemann surfaces. Then choosing the projective structure to be some canonical one, e.g. one
related to the Schottky or quasi-Fuchsian uniformization, one can get a functional on the (appropriate) moduli 
space itself. Since projective structures on Riemann surfaces are intimatly connected to
hyperbolic 3-manifolds, such Liouville functionals of projective structure turn out to
be related to the geometry of the corresponding 3-manifolds, and this is where the
notion of the renormalized volume enters into the story. 

However, the first examples of non-trivial Liouville functionals on the Teichm\"uller
space of Riemann surfaces had nothing to do with 3-manifolds, and worked entirely in
the 2-dimensional setting. Thus, the first construction for the moduli space of
surfaces of genus $g$ was proposed in \cite{TZ-Schottky}. The
idea was to base the construction on the so-called Schottky uniformization
of the surface $S$. The latter is unique given a {\it marking} of $S$, which
is a choice of $g$ generators of $\pi_1(S)$ corresponding too disjoint simple
closed curves on $S$. The surface is then cut along the
generators chosen, which results in a $2g$-holed sphere. Thus, in the Schottky uniformization,
the Riemann surface $S$ is represented as a quotient of the complex plane
under the action of the Schottky group $\Gamma$, which is a group freely generated
by $g$ (loxodromic, i.e. non-elliptic) elements of ${\rm PSL}(2,\C)$. A
fundamental region for the action of the Schottky group on $\C$ is
the exterior of a set of $2g$ Jordan curves that get mapped
into each other pairwise by the generators of $\Gamma$. The flat metric 
$|dz|^2$ of $\C$ can now be used as the reference structure for 
constructing a non-trivial Liouville functional, in place of $h_0$ above. This functional is
given by essentially the same quantity as in (\ref{liouv-1}), with the flat
metric Laplacian and $K_0=0$. However, the Liouville field $\phi$ now has
rather non-trivial transformation properties under the action of $\Gamma$
(so that the metric $e^\phi|dz|^2$ can be pulled back to the quotient Riemann surface).
As a result, the integral of the term $|\nabla\phi|^2$ depends on a choice
of the fundamental region. Thus, in order to define the functional in an
unambiguous way one has to correct the ``bulk'' term (\ref{liouv-1}) with
a rather non-trivial set of boundary terms, see \cite{TZ-Schottky} for details.
After this is done, one gets a well-defined functional of the conformal structure of $S$, 
of the marking of $S$ that was used to get the Schottky uniformization 
as well as the Liouville field $\phi$. One of the main properties of this
functional is that the 
variation with respect to $\phi$ gives rise precisely to the Liouville 
equation $\Delta\phi=e^\phi$, which can then be shown to have a unique solution 
with the required transformation properties. Evaluating the Liouville functional 
on this solution gives a non-trivial functional on the Schottky space, i.e. 
the moduli space of marked (by $g$ curves) Riemann surfaces. 
One of the key properties of this functional is that
its first variation (when the moduli are varied) gives rise to a certain
holomorphic quadratic differential on the complex plane, which measures
the difference between the Schottky and Fuchsian (used as reference) projective
structures. From this one finds that the extremal Liouville functional
is the K\"ahler potential for the so-called Weil-Petersson metric on 
the Schottky space, see \cite{TZ-Schottky} for details (and \cite{wolpert-review} for
many key properties of the Weil-Petersson metrics).

Even though this point of view is not at all developed in \cite{TZ-Schottky},
it is convenient to think of the functional constructed as a special case of a 
more general Liouville functional on a surface equipped with a projective structure, 
specialized to the case when the projective structure in question is the Schottky
one. This point of view suggests that there are other Liouville functionals
out there, namely when one chooses the projective structure to be different.
Indeed, a Liouville functional of the same type, but this time for the
so-called quasi-Fuchsian (or even more generally arbitrary Kleinian) projective structure 
was constructed in \cite{Takhtajan:2002cc}. Since quasi-Fuchsian projective structures
are parametrized by the product of two copies of Teichm\"uller space $\cT_g$, one gets a functional on
$\cT_g\times \cT_g$. Similar to the story in the Schotty case, this functional
turns out to be a K\"ahler potential for the Weil-Petersson metric on $\cT_g$,
as its $(1,1)$ derivative with respect to moduli of the first copy of $\cT_g$
turns out to be independent of the point in the second copy, see below for
more details. 

In \cite{takhtajan-zograf:spheres}, a 
Liouville functional of the same type was constructed for the moduli space 
$\cM_{(0,n)}^{(\alpha_1,\ldots,\alpha_n)}$ of Riemann surfaced with $n$ conical singularities of 
given angle deficits $2\pi \alpha_i, i=1,\ldots,n$. Moreover, a set of metrics 
on the space $\cM_{(0,n)}^{(\alpha_1,\ldots,\alpha_n)}$ was introduced, which for all $\alpha_i=1$ 
(angle deficits $2\pi$) coincides with the usual Weil-Petersson metric on $\cM_{(0,n)}$. It was moreover shown 
that all these metrics are K\"ahler with the K\"ahler potential given by the
Liouville functional. A choice of the projective structure is implicit in this
case, as the natural projective structure on the complex plane with $n$ marked
points is used. 

Below we shall see how these Liouville functionals are intimately related to the geometry
of the hyperbolic 3-manifolds which the corresponding projective structures define. We note
however that, although hyperbolic surfaces with cone singularities fit nicely in the 3D 
geometric context considered here (see section 8), 
the question of a geometrical (3-dimensional) interpretation of our last example -- the 
functional on $\cM_{(0,n)}^{(\alpha_1,\ldots,\alpha_n)}$ constructed in 
\cite{takhtajan-zograf:spheres} -- remains open.

\subsection{The AdS-CFT correspondence and holography\index{AdS-CFT}\index{holography}}

An initially unrelated development occurred in the context of
the so-called AdS/CFT correspondence\index{AdS/CFT correspondence}\index{correspondence!AdS/CFT} 
of string theory, see \cite{Witten:1998qj}, in which 
asymptotically-hyperbolic manifolds play the key role. These are non-compact,
infinite volume (typically Einstein, or Einstein with non-trivial ``stress-energy tensor''
on the right hand side of Einstein equations) Riemannian manifolds 
that have a conformal boundary, near which
the manifold looks like the hyperbolic space of the corresponding dimension. Most interesting
for physics is the case of five-dimensional asymptotically-hyperbolic manifolds,
for in this case the five-dimensional gravity theory induces (or, as physicists say,
is dual to) a certain non-trivial gauge theory on the boundary, see \cite{Witten:1998qj}.
However, the simplest situation is that in three dimensions, where the 
Einstein condition implies constant curvature and one is led to consider simply a
hyperbolic manifold $M$. 

Given such an asymptotically-hyperbolic manifold, 
physicists are interested in computing the Einstein-Hilbert functional 
(given by an integral of the scalar curvature plus a multiple of the volume form) 
on the metric of $M$. For the hyperbolic metric, the integrand
reduces to a multiple of the volume form, so one is computing the volume
of the hyperbolic manifold $M$, which is infinite. However, physicists
are masters of extracting a finite answer from a divergent expression.
And so it was found that in many situations there is a ``canonical'' way to
extract a finite answer by ``regularizing'' the divergent volume, 
then subtracting the divergent contribution, and finally taking a limit
to ``remove the regulator''. This procedure, somewhat mysterious for a mathematician,
was applied in the context of Schottky three-dimensional
hyperbolic manifolds (quotients of the hyperbolic space by the
action of a Schottky group) in \cite{Krasnov:2000zq}, with a somewhat unexpected
result. Namely, it was shown that the ``renormalized volume'' of a 
Schottky manifold, as a function of the conformal factor for a metric in the 
conformal class of the hyperbolic boundary, is equal to the Schottky Liouville functional
as defined in \cite{TZ-Schottky}. When one takes the conformal factor corresponding
to the canonical metric of curvature $-1$ one gets the Liouville functional on
the Schottky moduli space defined and studied in \cite{TZ-Schottky}. In view of the results of
this reference, one finds that the renormalized volume, which is a purely three-dimensional
quantity in its definition, equals the K\"ahler potential for the Weil-Petersson metric
on the Schottky space, which is a purely two-dimensional quantity. One thus
gets an example of what physicists like to refer to as a ``holographic''
correspondence (a relation between one quantity (or even theory) in $n+1$ dimensions
and another in $n$ dimensions). 

The methods of \cite{Krasnov:2000zq}, originally applicable only in the context of
(classical) Schottky manifolds, were generalized and applied to arbitrary
Schottky, quasi-Fuchsian and even Kleinian manifolds in \cite{Takhtajan:2002cc},
with the end result being always the same: the renormalized volume turns
out to be equal to the (appropriately defined) Liouville functional, and
the latter is shown to be the K\"ahler potential on the corresponding
moduli space. Even prior to work \cite{Takhtajan:2002cc}, the set of
ideas building upon the K\"ahler property of the Weil-Petersson metric on 
Teichm\"uller space (and closely related to the renormalized volume
ideas, as was later shown in \cite{Takhtajan:2002cc}) was used in 
\cite{McMullen} for a proof of the K\"ahler hyperbolic property of the moduli space.

The above story makes it clear that there is a deep relation between the
geometry of the Teichm\"uller space of a Riemann surface $S$ and the
geometry of hyperbolic three-manifolds that realize $S$ as its conformal
boundary. This relation was recently demonstrated from a more geometrical 
perspective in \cite{volume}, where it was shown that the
key property of the K\"ahler potential, namely, that its first variation is
given by a certain canonical quadratic differential on $S$, is a simple
consequence of the Schl\"afli-like formula of \cite{sem-era}. This geometrical
perspective on the renormalized volume (and the Liouville functional) 
also made it clear that such a quantity may be defined not only for
the considered in the literature Schottky, quasi-Fuchsian and 
Kleinian projective structures on the boundary, but, in fact, for
an arbitrary projective structure. This idea leads to the notion of
relative volume, defined and studied in \cite{cp}. It is
on this more geometrical and, we believe, very simple, perspective
on the renormalized volume that we would like to emphasize in this review.

\subsection{The definitions of the renormalized volume\index{volume!renormalized}}

Let $M$ be a convex co-compact hyperbolic manifold, for instance a quasifuchsian
manifold. The hyperbolic volume of $M$ is infinite, but an interesting finite quantity can
be extracted by a procedure that physicists refer to as {\it renormalization}. This
proceeds by introducing a ``regulator'' that makes the quantity of interest finite, then
removing the divergent contribution and then removing the regulator. For the case at
hand an appropriate regulator is given by equidistant surfaces. Thus, the
first step in the definition of the renormalized 
volume\index{renormalized volume}\index{volume!renormalized} of $M$ is
to define a volume depending on an equidistant foliation $F$ of $M$ in the complement
of a compact convex subset $N$. By an equidistant foliation we mean a foliation of 
$M\setminus N$ by closed, smooth, convex surfaces, so that, in each connected
component of $M\setminus N$, the surfaces are pairwise at constant distance.
Since the foliation is equidistant, it is uniquely determined by $N$.

Section 2 gives two different but equivalent definitions of a volume associated
to $N$ in $M$. The first definition is in terms of the asymptotic behavior
of the volume of the set of points at distance at most $\rho$ from $N$ as 
$\rho\rightarrow \infty$. This definition can be surprising but it extends
in higher dimensions in the setting of conformally compact, Poincar\'e-Einstein
metrics, see \cite{graham-witten} for an original reference and \cite{anderson:geometric} for a review. 

The other definition is simpler, but limited to $3$ dimensions. It is in terms of 
the volume of $N$, corrected by a term involving the integral of the mean 
curvature of the boundary of $N$, as
$$ W(M,N) = V(N)-\frac{1}{4} \int_{\dr N} H da~. $$

In 3 dimensions a convex co-compact hyperbolic manifold is completely specified by
the conformal structure of its boundary $\partial M$. It then turns out that the dependence of $W(N,M)$ 
on $N$ is just the dependence on a metric in the conformal class of the boundary. Moreover,
when this metric is varied the functional $W(N,M)$ reaches an extremum on metrics
of constant (negative) curvature. Thus, $W(N,M)$ is nothing but the Liouville
functional whose two-dimensional realizations have been discussed above. Note that
the hyperbolic 3-manifold $M$ in which the volume is computed comes equipped with
a projective structure on $\partial M$, and explains why the Liouville action in
all cases needed a projective structure to be defined. 

\subsection{The first and second fundamental forms at infinity}

There is a natural description of the connected components of the complement of 
$N$ in terms of the induced metric and second fundamental forms $I, \II$ of 
the corresponding boundary component of $N$. 
Section 3 contains an alternate description, in terms of naturally defined
``induced metric'' and ``second fundamental form'' $I^*,\II^*$ at infinity in the same 
end. There are simple transformation formulas from $I,\II$ to $I^*,\II^*$
and conversely. The conformal class of $I^*$ is the conformal class at infinity
of $M$. Moreover, $I^*$ and $\II^*$ satisfy the Codazzi equation and an analog
of the Gauss equation, which involves the trace of $\II^*$ instead of the
determinant of $\II$ as in the usual Gauss equation for surfaces in $H^3$.

A subtle point, explained at more length in Section 3, is that, as we have already
mentioned, not only $N$ determines the metric at infinity $I^*$, but $I^*$ also determines uniquely
$N$. In general a metric $I^*$ in the conformal class at infinity does not come
from a choice of some convex subset $N$ of $M$.  However, for any such metric
$I^*$, there is associated an equidistant foliation of a neighborhood of infinity
in $M$. This is sufficient to define the renormalized volume $W$, and this makes
$W$ a function of $M$ and $I^*$, rather than of $M$ and $N$ as defined above. 

The modified Gauss and Codazzi
equations for $I^*,\II^*$ are described in Section 3. They have an interesting
consequence. When $K^*$ is constant, the trace $H^*$ of $\II^*$ is also constant,
so that the traceless part $\II^*_0$ of $\II^*$ also satisfies
the Codazzi equation $d^{\nabla^*}\II^*_0=0$, where $\nabla^*$ is the
Levi-Civita connection of $I^*$. It then follows (by an argument
going back to Hopf \cite{hopf}) that $\II^*_0$ is the real part
of a holomorphic quadratic differential. In other words, $\II^*_0$
is canonically identified with a cotangent vector in $T^*_{[I^*]}\cT_{\dr M}$,
where $[I^*]$ is the conformal class of $I^*$, a property that is
going to play an important role below.

\subsection{Variations formulas}

The function $W(M,N)$ has a simple first-order variation formula in terms of
the data on the boundary of $N$, recalled in Section 4. For any first-order
deformation of $M$ or of $N$ in $M$,
$$ \delta W = \frac{1}{4}\int_{\dr N} \delta H + \langle\delta I,\II_0\rangle da, $$
where $H=\tr_I\II$ is the mean curvature of the boundary and $\II_0$ is the
traceless part of the second fundamental form. 

A lengthy but direct computation then shows that the same formula (except for the sign)
holds (when this function is interpreted as $W(M,I^*)$) in terms of the data at infinity:
\begin{equation} \label{eq:dW}
 \delta W = - \frac{1}{4}\int_{\dr N} \delta H^* + 
\langle\delta I^*,\II^*_0\rangle da,   
\end{equation}
where $H^*=\tr_{I^*}\II^*$ and $\II^*_0=\II^*-(H^*/2) I^*$. 
A version of this formula that is valid in any dimensions has been first obtained in
\cite{skenderis-solod} by a direct computation of the variation of the regularized
volume, see also \cite{anderson:geometric} for a more mathematical exposition. A geometrical
viewpoint adopted in this review (that originates in \cite{volume}) 
interprets this variational formula as a version of the Schl\"afli formula.

\subsection{Maximization}

A first consequence of Equation (\ref{eq:dW}), together with a simple
integration by parts argument which is explained in Subsection 4.3, is
that the only critical point of $W$ over the space of metrics of fixed
area in a conformal class on $\dr M$ is attained for the unique metric of constant 
curvature. Conversely, metrics of constant curvature are critical points
of the restriction of $W$ to metrics of fixed area in a conformal class,
and those critical points happen to be always non-degenerate minima.

This leads to the definition of the renormalized
volume of $M$ (no choice of $I^*$ involved) as the maximum of $W(M,I^*)$,
obtained precisely when the metric at infinity, $I^*$, has constant
curvature $K^*=-1$. We call $W(M)$ this function. By the Ahlfors-Bers
theorem \cite{ahlfors-bers}, $M$ is uniquely determined by its conformal class at
infinity, so $W$ can also be considered as a function from the 
Teichm\"uller space $\cT_{\dr M}$ of the boundary to $\R$.

The second step is to vary the conformal class on $\dr M$ while 
considering only constant curvature metrics. Then it follows
directly from (\ref{eq:dW}) that 
$$ dW(I^{*\prime}) = \frac{-1}{4}\int_{\dr M} \langle I^{*\prime},\II^*_0\rangle da^*~, $$
which means that $dW$ is identified (up to the factor $-1/4$)
with $\II_0^*$, considered as a 1-form on $\cT_{\dr M}$. 
 
\subsection{The renormalized volume and the Schwarzian 
derivative\index{Schwarzian derivative}\index{derivative!Schwarzian}}
\label{ssc:19}

There is another way to give a geometric meaning to $\II^*_0$. Given
a convex co-compact hyperbolic metric on $M$, it defines on the boundary
at infinity $\dr_\infty M$ a complex projective structure $\sigma$, see \cite{dumas-survey}. 
Let $c$ be the underlying complex structure, so that $c$ is the complex
structure at infinity of $M$. There is another special complex projective structure
on $\dr M$, obtained by applying the Riemann
uniformization theorem to $(\dr M, c)$, we call it $\sigma_0$. The 
image by the developing map of $\sigma_0$ of each connected component
of $\dr M$ is a disk, so $\sigma_0$ is called the Fuchsian complex
projective structure associated to $c$. 

Let $\phi$ be the map between $(\dr M,\sigma_0)$ and $(\dr M,\sigma)$. 
By construction $\phi$ is holomorphic, so that we can consider its 
{\it Schwarzian derivative} $\cS(\phi)$, which is a holomorphic quadratic
differential on $(\dr M,c)$. This holomorphic quadratic differential
can also be obtained as the difference between the projective structures
$\sigma$ and $\sigma_0$. We have:

\begin{prop}
$\II^*_0 = -Re(\cS(\phi))$.
\end{prop}

\subsection{The renormalized volume and Kleinian reciprocity\index{reciprocity!kleinian}}

Suppose now that $M$ is a quasifuchsian manifold, that is, it is convex co-compact and
homeomorphic
to the product of a closed oriented surface $S$ of genus at least $2$ with an interval.
Let $\cT$ be the Teichm\"uller space of $S$, and let $\cTb$ be the Teichm\"uller
space of $S$ with the opposite orientation. 
Given two complex structures $c_+\in \cT, c_-\in \cTb$, there is by Bers' Double
Uniformization Theorem a 
unique hyperbolic metric on $M$ such that the complex structure at infinity on the
upper boundary $\dr _+M$ of $M$ is $c_+$, while the complex structure at infinity on the lower
boundary $\dr_-M$ of $M$ is $c_-$. 

Given a quasi-Fuchsian manifold $M$, we can also consider the corresponding complex projective 
structure $\sigma_+$ on $\dr_+M$, and the Fuchsian complex projective structure 
$\sigma_{0,+}$ on $\dr_+M$ obtained by applying the Riemann uniformization theorem to the complex 
structure $c_+$. This defines, through the Schwarzian derivative construction recalled above, a
holomorphic quadratic differential $q_+$ on $(\dr_+M,c_+)$, and the real part 
of $q_+$ defines a cotangent vector $\beta_+(c_-,c_+)\in T^*_{c_+}\cT$. The 
same construction for $\dr_-M$ yields a cotangent vector $\beta_-(c_-,c_+)\in
T^*_{c_-}\cTb$.  

McMullen's quasifuchsian reciprocity\index{quasifuchsian reciprocity}\index{reciprocity!quasifuchsian} 
\cite{McMullen} gives a subtle relation between the ways
the complex projective structures on the two boundary components
vary when the complex structure at infinity changes. In the following
statement we consider a fixed conformal structure $c_+$ and vary $c_-$ 
(resp. fix $c_-$ and vary $c_+$), $\beta_+$ is then a map from 
$\cTb$ to $T^*_{c_+}\cT$ (resp. $\beta_-$ is a map from 
$\cT$ to $T^*_{c_-}\cTb$). We call $D\beta_+$ (resp. $D\beta_-$) the
differential of this map.

\begin{thm}[McMullen \cite{McMullen}] \label{tm:mcm}
The tangent maps 
$$ D\beta_+(\cdot, c_+):T_{c_-}\cTb\rightarrow T^*_{c_+}\cT $$
and 
$$ D\beta_-(c_-, \cdot):T_{c_+}\cT\rightarrow T^*_{c_-}\cTb $$
are adjoint.
\end{thm}

This can be stated in different terms using the standard cotangent bundle symplectic structure on 
$T^*\cT_{\dr M}$, which we will call $\omega_*$ here. Given $(c_-,c_+)\in \cTb\times
\cT$, we call $\beta(c_-,c_+)=(\beta_-(c_-,c_+,),\beta_+(c_-,c_+))\in T^*_{c_-}\cTb\times
T^*_{c_+}\cT$, so that $\beta(c_-,c_+)\in T^*_{(c_-,c_+)}\cT_{\dr M}$. Thus $\beta$
defines a section of $T^*\cT_{\dr M}$.
Theorem \ref{tm:mcm} is a direct consequence of the following simpler statement.

\begin{thm} \label{tm:mcm2}
The image $\beta(\cT_{\dr M})$ is Lagrangian in $(T^*\cT_{\dr M},\omega_*)$.
\end{thm}

In this form, the statement extends as it is to a much more general setting
of convex co-compact (or geometrically finite) hyperbolic 3-manifolds. 

The proof of Theorem \ref{tm:mcm} from Theorem \ref{tm:mcm2} is straightforward,
as is the proof of Theorem \ref{tm:mcm2} from the first-order variation of the 
renormalized volume as described above. Both are done in Subsection \ref{ssc:renormalized}.
It is proved there that Theorem \ref{tm:mcm2} is actually equivalent to Theorem \ref{tm:mcm} 
along with the following proposition, which uses the same notations as Theorem \ref{tm:mcm}
and is another direct consequence of Theorem \ref{tm:mcm2}.

\begin{cor} \label{cr:closed}
For fixed $c_-$, $\beta_{+}(c_-,\cdot)$ is a closed 1-form on $\cT$. 
\end{cor}

Part \ref{ssc:ccore} describes a second proof of Theorem \ref{tm:mcm2}, based on the
geometry of the convex core and on the fact that the grafting map is symplectic
(Theorem \ref{tm:cp} below). The idea here is to prove first that the data on the
boundary of the convex core defines, as the representation varies, a Lagrangian 
submanifold of the cotangent bundle of Teichm\"uller space, understood here in
terms of hyperbolic metrics and measured laminations. But the natural map 
sending this data on the boundary of the convex core to the corresponding 
data at infinity is symplectic, and Theorem \ref{tm:mcm2} follows.

In Part \ref{ssc:kerckhoff} we describe a third argument, due to Kerckhoff, 
which works in the setting of deformations of the holonomy representation of
the fundamental group of $M$. This uses topological arguments, more precisely
a long exact sequence and Poincar\'e duality between cohomology spaces with
value in an $sl(2,\R)$-bundle over $M$. The equivalence with Theorem \ref{tm:mcm2}
is made through a result of Kawai \cite{kawai} connecting the Goldman symplectic
form on the space of complex projective structures on a surface to the cotangent
symplectic structure on the cotangent bundle of Teichm\"uller space.

\subsection{The renormalized volume as a K\"ahler 
potential\index{potential!K\"ahler}\index{K\"ahler potential}}

A rather direct consequence of McMullen's quasifuchsian reciprocity, explained in
Section \ref{sc:kaehler}, is that although $\beta_+(c_-,c_+)$ clearly depends on 
$c_-$, its exterior (anti-holomorphic, because $\beta_+(c_-,c_+)$ is a holomorphic one-form
on $\cT$) differential does not depend on $c_-$. So this exterior
differential can be computed explicitly in the simplest possible case -- when 
$c_-=c_+$, that is, in the neighborhood of  a ``Fuchsian'' hyperbolic manifold.
This leads to a key result on the renormalized volume, namely that it is
a K\"ahler potential for the Weil-Petersson metric on Teichm\"uller space.

\subsection{The relative volume of hyperbolic ends and the grafting map\index{volume!relative}\index{grafting}}

So far we have considered the renormalized volume in the context of hyperbolic convex
co-compact 3-manifolds. Such manifolds come equipped with a complex projective
structure on each boundary component, and the renormalized volume we have
discussed can be said to be the Louville functional for the corresponding
projective structure, as discussed in the beginning of this section.
Section \ref{sc:relative} is centered on a notion of renormalized volume,
or Liouville functional, that is defined for an arbitrary projective structure
on a Riemann surface $S$. This has been developed in \cite{cp}. The main idea is to use a variant
of the renormalized volume, called the {\it relative volume} of a hyperbolic
end, and then use an analog of Theorem \ref{tm:mcm2} to obtain results on the grafting map. 

Thus, consider a closed surface
$S$ of genus at least $2$. We denote by $\ML_S$ the space of measured 
geodesic laminations on $S$, see e.g. \cite{bonahon-geodesic}, and by $\CP_S$ the space of
complex projective structures (or $\C P^1$-structures) on $S$, see e.g. \cite{dumas-survey}.

The grafting map $Gr:\cT_S\times \ML_S\rightarrow \CP_S$ was defined by Thurston,
see e.g. \cite{dumas-wolf,dumas-survey}. It can be described (superficially) 
as follows. Let $m\in \cT_S$ be a hyperbolic
metric, and let $l\in \ML_S$ be a measured lamination with support a disjoint union of
closed curves $c_1,\cdots, c_n$, each with a positive weigth $w_1, \cdots, w_n$. 
The ``grafted metric'' on $S$ is obtained by cutting open $(S,m)$ along each of 
the $c_i$ and gluing in a flat strip of width equal to $w_i$. Then $Gr(m,l)$ is 
the complex projective structure naturally associated to this metric. The map
$Gr$ extends by continuity to a map defined on all measured laminations (not only
those supported by multicurves). 

Thurston showed that the grafting map is a homeomorphism. His
proof, although quite subtle, relies on a simple geometric idea, which is simpler
to explain for a complex projective structure $\sigma \in \CP_S$ whose developing
map $\dev:\St\rightarrow \C P^1$ is injective. In this case, the boundary of the
convex hull in $H^3$ of $\C P^1\setminus \dev(\St)$ -- here $\C P^1$ is identified
with the boundary at infinity of $H^3$ -- is a convex pleated surface, on which
$\pi_1(S)$ acts properly discontinuously. The quotient surface is endowed with a 
hyperbolic metric (its induced metric) and a measured lamination (its pleating
lamination). This gives an element $(m,l)\in \cT_S\times \ML_S$, which is the
inverse image of $\sigma$ under the grafting map. 

In this picture -- which can be extended to the case where the developing map
of $\sigma$ is not injective -- $\pi_1(S)$ acts properly discontinously on
the connected component bounded at infinity by $\dev(\St)$
of the complement in $H^3$ of the convex hull of 
$\C P^1\setminus \dev(\St)$, and the quotient is a {\it hyperbolic end}. 
In Section \ref{sc:relative} we explain how to define the {\it relative 
volume} of such a hyperbolic end, and show that it satisfies a simple
first-order variation formula, involving both a term ``at infinity'' similar
to the one which we already mentioned for the renormalized volume, and
a  term on the ``compact'' boundary, involving the hyperbolic metric and
the measured pleating lamination, which is very close to a Schl\"afli-type
formula for the convex core of a quasifuchsian manifold, discovered by
Bonahon \cite{bonahon-variations}. 

There is a natural identification between $\cT_S\times \ML_S$ and the
cotangent space $T^*\cT_S$, obtained by considering the differential of the
length of a measured lamination as a cotangent vector. Using this map, 
the grafting map can be considered as a map from $\Grb:T^*\cT_S\rightarrow 
\CP_S$, and both sides are naturally symplectic manifolds ($\CP_S$ is 
actually a complex symplectic manifold, but we consider here only
the real part of its complex symplectic structure). It is proved in 
\cite{cp} -- using mostly tools from Bonahon's work \cite{bonahon,bonahon-variations} 
-- that this map is $C^1$ smooth (it is however not $C^2$). 

\begin{thm} \label{tm:cp}
$(1/2)\Grb$ is symplectic. 
\end{thm}

The proof is a direct consequence of the first-order variation formula for
the relative volume of hyperbolic ends, similar to the proof of Theorem 
\ref{tm:mcm} that follows from the first-order variations of the renormalized
volume. 

\subsection*{Acknowledgements}

We are grateful to Steve Kerckhoff for interesting conversations and for
allowing us to present here the content of part \ref{ssc:kerckhoff}. We 
would also like to thank Brice Loustau for pointing out an error in a previous
version of this text, and for stimulating conversations on the questions 
considered here.

\section{Two definitions of the renormalized volume}

\subsection{The renormalized volume}

As mentioned in the introduction, the renormalized volume was motivated by
physical considerations. This notion was first envisaged and used in the more general 
context of conformally compact (in particular Einstein) manifolds in an arbitrary number of
dimensions, and only later specialized to the three-dimensional setting considered here. In the 
more general case its definition uses a foliation of a particular type close to infinity
(associated to a ``defining function'' of the boundary). The "canonical" renormalized
volume independent of which foliation is used is then either the constant term (in even dimensions) or 
the logarithmic term (in odd dimensions) in the asymptotic expansion of the volume in terms of
the parameter of the foliation. In this review, however, we are interested in the constant
term in this asymptotic expansion for the odd-dimensional (3d) manifolds. This quantity
depends on the choice of a metric in the conformal class at infinity. The higher-dimensional applications can be 
found e.g. in \cite{graham-witten,herzlich-renormalized,anderson:geometric}, see also 
\cite{skenderis,skenderis-solod}, while this review concentrates on a simple(r), but 
still extremelly rich case of three-dimensional manifolds.

\subsection{The renormalized volume via equidistant 
foliations\index{equidistant foliation}\index{foliation!equidistant}}

The limiting procedure via which the
volume is defined can be somewhat de-mystified by considering for
regularization a family of surfaces equidistant to a given one, following an
idea already used by C. Epstein \cite{epstein-duke} (and more recently 
put to use in \cite{minsurf}). Thus, the main idea
is to obtain the renormalized volume by taking
a convex domain $N\subset M$, and compute the renormalized volume of $M$ {\it
  with respect to $N$} as 
\be\label{rv}
V_R(M,N)=V(N) + \lim_{\rho\to\infty} \left( V(\partial N, \partial N_\rho) - 
(1/2) A(\partial N_\rho) - \sum_i 2\pi \rho(g_i-1) \right)
\ee
where $V(\partial N, \partial N_\rho)$ is the volume between the boundary $\partial N$
of the domain $N$ and the surface $\partial N_\rho$ located a distance $\rho$ from $\partial N$.
The quantity $A(\partial N_\rho)$ is the area of the surface $\partial N_\rho$, the sum
in the last term is taken over all boundary components of $M$ and $g_i$ are the genera of
these boundary components. The convexity of
the domain $N$ ensures that the equidistant surfaces $\partial N_\rho$ exist all the way to
infinity. This ensures that the limit $\rho\to\infty$ can be taken. This 
limit exists and can be computed in terms of the volume of $N$, 
see below for the corresponding expressions. The limiting procedure 
used in \cite{Krasnov:2000zq}, \cite{Takhtajan:2002cc} is an example
of the limiting procedure described above, for the Epstein surfaces \cite{Eps} used in these
references are equidistant. 

\subsection{Renormalized volume as the $W$-volume}

The following formula for the renormalized volume 
(\ref{rv}) can be shown by an explicit (simple) computation:
\be
V_R(M, N) = W(N) - \sum_i \pi(g_i-1)~, 
\ee
where the sum in the last term is taken over all the boundary components. Here
the W-volume is defined as
\be\label{Wvol}
W(N) := V(N) - \frac{1}{4} \int_{\dr N} H da.
\ee
This formula for $V_R$ is a special case of a formula found by C. Epstein 
\cite{epstein-duke} for the renormalized volume of hyperbolic manifolds 
in any dimension. 

Thus, the renormalized volume of $M$ with respect to $N$ is, apart from an
uninteresting term given by a multiple of the Euler characteristic of the
boundary, just the W-volume of the domain $N$. It then makes sense to study this
geometrical $W$-volume instead. Note already that $W(N)$ is not equal to the 
Hilbert-Einstein functional of $N$ with its usual boundary term, it
differs from it in the coefficient of the boundary term.

\subsection{Self-duality}

One of the interesting properties of the W-volume is that it
is self-dual. Thus, we recall that the Einstein-Hilbert functional
\be\label{gr-2}
I_{EH}(N) := V(N) - \frac{1}{2} \int_{\dr N} H da
\ee
for a compact domain $N\subset H^3$ of hyperbolic space
(note a different numerical factor in front of the second term) is nothing but the
dual volume. Indeed, we recall that there is a duality between objects in $H^3$ and
objects in $dS_3$, the 2+1 dimensional de Sitter space. Under this duality
geodesic planes in $H^3$ are dual to points in $dS_3$, etc. This duality between
domains in the two spaces is easiest to visualize for convex polyhedra
(see \cite{RH}), but the duality works for general domains as well. The fact 
that (\ref{gr-2}) is the volume of the dual domain 
is a simple consequence of the Schl\"afli formula, see below.

Thus, we can write:
\be
{}^*V(N) = V({}^*N) = V(N) - \frac{1}{2} \int_{\dr N} H da
\ee
for the volume of the dual domain. This immediately shows that
\be
W(N) = V(N) - \frac{1}{4} \int_{\dr N} H da = \frac{V(N)+{}^*V(N)}{2}.
\ee
Thus, the $W$-volume is self-dual in that this quantity for $N$ is equal to this quantity
for the dual domain ${}^*N$: $W(N)=W({}^*N)$. 

\subsection{The $W$-volume and the Chern-Simons 
formulation\index{Chern-Simons formulation}\index{formulation!Chern-Simons}} 

An interesting remark is that
there is a very simple expression for the $W$-volume in terms of the so-called
Chern-Simons formulation of 2+1 gravity \cite{Witten:2+1}. Let us briefly review this formulation.
In the so-called first order formalism for gravity the independent variables are not components
of the spacetime metric but instead the triads and the spin connection. Thus, for
Riemannian signature gravity in 3 spacetime dimensions let us introduce a 
collection of 3 one-forms $\theta^i, i=1,2,3$ such that the spacetime interval can be written
as $ds^2= \sum_i \theta^i\otimes \theta^i$. We can now construct from $\theta^i$ an 
$\mathfrak{su}(2)$ Lie algebra-valued one-form by taking $\theta=\sum_i {\rm i} \theta^i \sigma^i$,
where the $\sigma^i$ are the $2\times 2$ anti-symmetric Pauli matrices 
$\sigma^i \sigma^j = \delta^{ij}{\rm Id} + {\rm i} \epsilon^{ijk}\sigma^k$. Using the
Lie algebra-valued form $\theta$ one can write the metric as $ds^2 = - (1/2){\rm Tr}(\theta\otimes\theta)$.

The Einstein-Hilbert action as a functional of the metric $g$ is given by 
\be\label{EH}
I_{EH}[g]=-\frac{1}{4} \int_M dv\, (R+2) - \frac{1}{2}\int_{\partial M} da \, H,
\ee
where $R$ is the Ricci scalar of $g$, $H$ is the mean curvature of the boundary, and $dv, da$ are the 
volume and area forms on $M$ and $\partial M$ respectively. When evaluated on a constant curvature
metric with $R=-6$ the Einstein-Hilbert action reduces to (\ref{gr-2}). The functional 
(\ref{EH}) can be re-written in a very simple form in terms of $\theta$ by introducing a 
spin connection $\omega$, which is locally an $\mathfrak{su}(2)$-valued one-form. The action is then:
\be
I_{EH}[\theta,\omega]=
\frac{1}{2}\int_M {\rm Tr}(\theta\wedge f(\omega) - \frac{1}{12} \theta\wedge \theta\wedge \theta)+
\frac{1}{2} \int_{\partial M} {\rm Tr}(\theta\wedge\omega).
\ee
Here $f(w)=d\omega+\omega\wedge\omega$ is the curvature of the spin connection $\omega$. When
one varies this action with respect to $\omega$ one obtains the equation $d_\omega \theta=0$,
where $d_\omega$ is the covariant derivative with respect to the connection $\omega$.
This equation can be solved for $\omega$ in terms of the derivatives of $\theta$. Once
one substitutes the solution back into the action one gets (\ref{EH}).

In contrast, the combination (\ref{Wvol}) that plays the role of the renormalized volume 
is obtained by evaluating on the constant curvature metric the following action:
\be
I_{W}[g]=-\frac{1}{4} \int_M dv\, (R+2) - \frac{1}{4}\int_{\partial M} da \, H.
\ee
This can be written in terms of the tetrad and the spin connection forms as follows:
\be\label{W-act}
I_{W}[\theta,\omega]=
\frac{1}{2}\int_M {\rm Tr}(\theta\wedge f(\omega) - \frac{1}{12} \theta\wedge \theta\wedge \theta)+
\frac{1}{4} \int_{\partial M} {\rm Tr}(\theta\wedge\omega).
\ee

One then notes that this is precisely the combination that appears in the Chern-Simons formulation.
Indeed, let us introduce the Chern-Simons action of an ${\mathfrak su}_{\C}(2)$ connection $A$
via:
\be
I_{CS}[A]= \frac{1}{4{\rm i}} \int_M {\rm Tr}(A\wedge dA + \frac{2}{3}A\wedge A\wedge A).
\ee
Now, defining the two ${\mathfrak su}_{\C}(2)$ connections:
\be
A = \omega+ \frac{i}{2}\theta, \qquad \bar{A} = \omega - \frac{i}{2}\theta
\ee
it is not hard to see that (\ref{W-act}) is given by:
\be
I_{W} = I_{CS}[A] - I_{CS}[\bar{A}],
\ee
with precisely the right boundary term that comes from having to integrate by parts.
In contrast, to obtain via the Chern-Simons formulation the combination (\ref{gr-2}) one
has to add a separate boundary term that is constructed from both $A,\bar{A}$. For
more details on the argument presented the reader is referred to \cite{Krasnov:2003}, 
see formula (3.7) of this reference as well as the related discussion. It would be of 
interest to understand the relation, if any, between the self-duality of the $W$-volume
and the fact that it has such a simple expression in the Chern-Simons
formulation.

\section{Description ``from infinity''}

\subsection{The metric at infinity\index{metric at infinity}}

In this section we switch from a description of the renormalized volume from
the boundary of a convex subset to the boundary at infinity of $M$. This
description from infinity is remarkably similar to the previous one from the
boundary of a convex subset.

\begin{lemma} \label{lm:asympt}
Let $M$ be a convex co-compact hyperbolic 3-manifold, and let $N\subset M$ be
compact and ``strongly'' convex with smooth boundary. Let $S_\rho$ be the
equidistant surfaces from $\dr N$. The induced metric on $S_\rho$ 
is asymptotic, as $\rho\rightarrow \infty$, to $(1/2) e^{2\rho} I^*$, where
$I^*=(1/2)(I+2\II+\III)$ is defined on $\dr N$.
\begin{proof}Follows from the following Lemma.
\end{proof}
\end{lemma}

\begin{lemma} \label{lm:base}
Let $S$ be a surface in $H^3$, with bounded principal curvatures, 
and let $I, B$ be the first fundamental 
form and the shape operator of $S$ correspondingly. Let $S_\rho$ be the
surface at  
distance $\rho$ from $S$. Then, for sufficiently small $\rho$ the induced
metric on $S_\rho$ is: 
\be \label{metric} I_\rho(x,y) =
I\left(\left(\cosh(\rho)E+\sinh(\rho)B\right)x,\left(\cosh(\rho)E+\sinh(\rho)B\right)y\right)~. 
\ee
Here $E$ is the identity operator.  
\end{lemma}
Note that this lemma also holds for a surface $S$ in any hyperbolic 3-manifold
$M$, not necessarily 
$H^3$. We also note that when the surface $S$ is convex, then the expression
(\ref{metric}) gives 
the induced metric on any surface $\rho>0$, where $\rho$ increases in the
convex direction. A proof of this lemma can be found in, e.g., \cite{minsurf}.

It is the metric $I^*$ that will play such a central role in what follows, so
we would like to state some of its properties.

\begin{lemma}
The curvature of $I^*$ is 
\be\label{curv*}
K^*:=\frac{2K}{1+H+K_e}~. 
\ee 
\end{lemma}

\begin{proof}
The Levi-Civita connection of $I^*$ is given, in terms of the Levi-Civita
connection $\nabla$ of $I$, by:
$$ \nabla^*_xy = (E+B)^{-1} \nabla_x((E+B)y)~. $$
This follows from checking the 3 points in the definition of the Levi-Civita 
connection of a metric:
\begin{itemize}
\item $\nabla^*$ is a connection.
\item $\nabla^*$ is compatible with $I^*$.
\item it is torsion-free (this follows from the fact that $E+B$ 
verifies the Codazzi equation: $(\nabla_x (E+B))y=(\nabla_y(E+B))x$). 
\end{itemize}

Let $(e_1, e_2)$ be an orthonormal moving frame on $S$ for $I$, and let
$\beta$ be its connection 1-form, i.e.:
$$ \nabla_xe_1 = \beta(x)e_2, ~\nabla_xe_2 = -\beta(x)e_1~. $$
Then the curvature of $I$ is defined as: $d\beta = -K da$. 

Now let $(e^*_1,e^*_2) := \sqrt{2}((E+B)^{-1}e_1,(E+B)^{-1}e_2)$; clearly it
is an orthonormal moving frame for $I^*$. Moreover the above expression of
$\nabla^*$ shows that its connection 1-form is also $\beta$.
It follows that $Kda = -d\beta = K^*da^*$, so that:
$$ K^* = K \frac{da}{da^*} = \frac{K}{(1/2)\det(E+B)} = \frac{2K}{1+H+K_e}~. $$
\end{proof}

We note that the metric $I^*$ is defined for any surface $S\subset
M$. However, it might have singularities (even when the surface $S$ is smooth)
unless $S$ is strictly {\it horospherically convex}, i.e., 
its principal curvatures are less than $1$ (which implies that it 
remains on the concave
side of the tangent horosphere at each point). If $S$ is a strictly 
horospherically convex surface $S$ embedded in
a hyperbolic end of $M$ then the metric $I^*$ 
is guaranteed to be in the conformal class of the  
(conformal) boundary at infinity of $M$. For
a general surface $S$ the ``asymptotic'' metric is not directly related to the
conformal infinity, and in particular, it does not have to be in the conformal
class of the boundary. 

\subsection{Second fundamental form at infinity\index{Second fundamental form at infinity}}

We have already defined the metric ``at infinity''. Let us now add to this
a definition of what can be called the second fundamental form at infinity.

\begin{df} Given a surface $S$ with first, second and third
fundamental forms $I,\II$ and $\III$, we 
define the first and second fundamental forms ``at infinity'' as:
\be\label{ff*}
I^* = \frac{1}{2}(I + 2\II + \III)= \frac{1}{2}(I+\II)I^{-1}(I+\II) =
\frac{1}{2}I((E+B)\cdot, (E+B)\cdot)~, \\
\nonumber 
\II^* = \frac{1}{2}(I-\III)=\frac{1}{2}(I+\II)I^{-1}(I-\II) = 
\frac{1}{2}I((E+B)\cdot, (E-B)\cdot)~.
\ee
\end{df}

It is then natural to define:
\be B^*:=(I^*)^{-1}\II^* = (E+B)^{-1}(E-B)~, \label{eq:def-B*} \ee
and 
$$ \III^*:=I^*(B^*\cdot, B^*\cdot) = I((E-B)\cdot, (E-B)\cdot)~. $$
An interesting point is that $I^*, \II^*$ and $\III^*$ determine the
full asymptotic development of the metric close to infinity: the 
induced metrics $I_\rho$ on the surfaces $S_\rho$ are given by Equation
(\ref{metric-infinity}) below. 
This extends Lemma \ref{lm:asympt}, and can be considered as an
analog of Equation (\ref{metric}).

Note that, for a surface which has principal curvatures strictly 
bounded between $-1$ and $1$, $\III^*$ is also a smooth metric and its
conformal class corresponds to that on the other component of the boundary at
infinity. This is a simple consequence of Lemma \ref{lm:base} and
the fact that when the principal curvatures are strictly bounded
between $-1,1$ the foliation by surfaces equidistant to $S$ 
extends all the way through the manifold $M$. Such manifolds were
called {\it almost-Fuchsian} in our work \cite{minsurf}.

As before, these definitions make sense for any surface, but it is only for a
convex surface (or more generally for a horospherically convex surface)
that the fundamental forms so introduced are guaranteed to have
something to do with the actual conformal infinity of the space.

\subsection{The Gauss and Codazzi equations at infinity\index{Gauss!equation}\index{Codazzi!equation}}

We also define $H^*:=\tr(B^*)$. The Gauss equation for ``usual'' surfaces in
$H^3$ is replaced by a slightly twisted version.

\begin{lemma} \label{rk:H*K*}
$H^*=-K^*$: the mean curvature at infinity is equal to minus 
the curvature of $I^*$.
\end{lemma}

\begin{proof}
By definition, $H^* = \tr((E+B)^{-1}(E-B))$. An elementary computation (for
instance based on the eigenvalues of $B$) shows that
$$ H^* = \frac{2-2\det(B)}{1+\tr(B)+\det(B)}~. $$
But we have seen (as Equation (\ref{curv*})) that $K^*=2K/(1+H+K_e)$.
The result follows because, by the Gauss equation, $K=-1+\det(B)$. 
\end{proof}

However, the ``usual'' Codazzi equation holds at infinity.

\begin{lemma} \label{rk:codazzi*}
$d^{\nabla^*}B^*=0$. 
\end{lemma}

\begin{proof}
Let $u,v$ be vector fields on $\dr_\infty M$. Then it follows from the
expression of $\nabla^*$ found above that:
\begin{eqnarray*}
(d^{\nabla^*}B^*)(x,y) & = & \nabla^*_x(B^*y) - \nabla^*_y(B^*x) - B^*[x,y] \\
& = & (E+B)^{-1}\nabla_x((E+B)B^*y) - \\
& & - (E+B)^{-1}\nabla_y((E+B)B^*x) - B^*[x,y]
\\ 
& = & (E+B)^{-1}\nabla_x((E-B)y) - \\
& & - (E+B)^{-1}\nabla_y((E-B)x) -
(E+B)^{-1}(E-B)[x,y] \\
& = & (E+B)^{-1}(d^\nabla (E-B))(x,y) \\
& = & 0~.
\end{eqnarray*}
\end{proof}

\subsection{Inverse transformations}

The transformation $I,\II \to I^*, \II^*$ is invertible. The inverse is
given explicitly, and the inversion formula exhibits a remarkable symmetry.

\begin{lemma} Given $I^*, \II^*$, the fundamental forms $I, \II$ such that
(\ref{ff*}) holds are obtained as:
\be\label{metric-rel}
I = \frac{1}{2} (I^*+\II^*)(I^*)^{-1}(I^*+\II^*) =
\frac{1}{2}I^*((E+B^*)\cdot, (E+B^*)\cdot))~, \\ \nonumber
\II = \frac{1}{2} (I^*+\II^*)(I^*)^{-1}(I^*-\II^*)=
\frac{1}{2}I^*((E+B^*)\cdot, (E-B^*)\cdot))~.
\ee
Moreover,
\be B = (E+B^*)^{-1}(E-B^*)~. \label{eq:B-inverse} \ee
\end{lemma}

Having an expression for the fundamental forms of a surface in terms of the
ones at infinity, one can re-write the metric of Lemma \ref{lm:base} induced 
on surfaces equidistant to $S$ in terms of $I^*,\II^*$.

\begin{lemma} The metric (\ref{metric}) induced on the surfaces equidistant to
  $S$  
can be re-written in terms of the fundamental forms ``at infinity'' as:
\be \label{metric-infinity}
I_\rho = \frac{1}{2} e^{2\rho} I^* + \II^* + \frac{1}{2} e^{-2\rho} \III^*~.
\ee
\end{lemma}
This lemma shows the significance of $\II^*$ as being the constant term of the metric.
This lemma also shows clearly that when the equidistant foliation extends all
the way through $M$ (i.e. when the principal curvatures on $S$ are in $(-1,1)$),
the conformal structure at the second boundary component of $M$ is that
of $\III^*=\II^* (I^*)^{-1} \II^*$. Thus, in this particular case of
almost-Fuchsian manifolds, the knowledge of $I^*$ on both boundary
components of $M$ is equivalent to the knowledge of $I^*, \II^*$ near either
component. In other words, $\II^*$ is determined by $I^*$. This statement
is more general and works for manifolds other than almost-Fuchsian.

\subsection{Fundamental Theorem of surface theory ``from infinity''}

Let us now recall that the Fundamental Theorem of 
surface theory states that
given $I,\II$ on $S$ satisfying the Gauss and Codazzi equations, 
there is a unique embedding of $S$ into the
hyperbolic space with induced metric and second fundamental form equal to 
$I$ and to $\II$. Then (\ref{metric}) gives an expression for the
metric on equidistant surfaces to $S$, and thus describes a hyperbolic manifold
$M$ in which $S$ is embedded, in some neighborhood of $S$. It would be 
possible to state a similar result for hyperbolic ends, uniquely 
determined by $I^*$ and $\II^*$ at infinity. But there is also an
analogous theorem, based on a classical result of Bers \cite{Bers},
in which the first (and only the first) form at infinity
is used. This can be compared with arguments used in \cite{horo}, where
a corresponding second fundamental form and the Gauss and Codazzi equations
``at infinity'' were introduced.

\begin{thm} 
Given a convex co-compact 3-manifold $M$, and 
a metric $I^*$ (on all the boundary components of $M$) in the conformal class
of the boundary,  
there is a unique foliation of each end of $M$ by convex equidistant surfaces
$S_\rho\subset M$ such that 
$(1/2)(I_\rho + 2\II_\rho + \III_\rho)=e^{2\rho} I^*$, 
where $I_\rho, \II_\rho, \III_\rho$ are the fundamental forms of $S_\rho$.
\end{thm}

\begin{remark} {\rm Note that one does not need to specify $\II^*$. The first
fundamental form $I^*$ (but on all the boundary components) is sufficient.}
\end{remark}

\begin{proof}
The surfaces in question can be
given explicitly as an embedding of the universal cover $\tilde{S}$ of $S$
into the hyperbolic space.  
Thus, let $(\xi,y), \xi>0, y\in \C$ be the usual upper half-space model
coordinates of $H^3$. Let us write the metric at infinity as 
\be\label{I*-Liouv}
I^* = e^{\phi} |dz|^2,
\ee
where $\phi$ is the Liouville
field covariant under the action of the Kleinian group giving $M$ on $S^2$. 
The surfaces
are given by the following set of maps: $\Eps_\rho: S^2 \to H^3, z\mapsto (\xi,y)$
(here $\Eps$ stands for Epstein, 
who described these surfaces in \cite{Eps}):
\be 
\xi = \frac{\sqrt{2} e^{-\rho} e^{-\phi/2}}{1+(1/2)e^{-2\rho} e^{-\phi}
  |\phi_z|^2}, \\ \nonumber 
y = z + \phi_{\bar{z}} \frac{e^{-2\rho} e^{-\phi}}{1+(1/2)e^{-2\rho} e^{-\phi}
  |\phi_z|^2}. 
\ee
As is shown by an explicit computation, the metric induced on the 
surfaces $S_\rho$ is given by
(\ref{metric-infinity}) with 
\be
\II^* = \frac{1}{2}(\theta dz^2 + \bar{\theta} d\bar{z}^2) + \phi_{z\bar{z}} dz d\bar{z}, \\ 
\label{theta}
\theta = \phi_{zz} - \frac{1}{2} (\phi_z)^2.
\ee
Thus, we see that $\II^*$ is determined by the conformal factor in (\ref{I*-Liouv}).
\end{proof}

\begin{remark} {\rm This theorem implies that the renormalized volume only depends on $I^*$. Indeed,
the foliation $(S_\rho)$ of the ends does depend only on $I^*$, and this foliation can be
used for regularization and subtraction procedure. Then the fact that the $W$-volume is
essentially the renormalized volume implies that the $W$-volume is a functional of $I^*$ only.
In the next section we will find a formula for the first variation of this functional.}
\end{remark}

\begin{cor} If the principal curvatures at infinity (eigenvalues of $B^*$) are
  positive the map $\Eps_\rho$ is a homeomorphism onto its image for any $\rho$.
\end{cor}

\begin{proof}
We first note that the map $\Eps_\rho$ is not always a 
homeomorphism, and the surfaces $S_\rho$ are not
necessarily convex, but for sufficiently large $\rho$ 
both things are true. A condition that guarantees that 
$\Eps_\rho$ is a homeomorphism for any $\rho$ is stated 
above. This condition can be
obtained from the requirement that the principal 
curvatures of surfaces $S_\rho$ are in $[-1,1]$. Let us
consider the surface $S:=S_{\rho=0}$ the first and 
second fundamental forms of which are given by 
(\ref{metric-rel}) (this immediately follows from (\ref{metric-infinity})). 
The shape operator of
this surface is then given by $B=(E+B^*)^{-1}(E-B^*)$. It is then 
clear that the principal
curvatures of $S$ are given by $k_i = (1-k_i^*)/(1+k_i^*)$, where the
$k_i^*$ are the ``principal curvatures''
(eigenvalues) of $B^*$. The latter are easily shown to be given by
\be
k^*_{1,2} = e^{-\phi}\left(\phi_{z\bar{z}}\pm \sqrt{\theta\bar{\theta}}\right).
\ee
It is now easy to see that the condition $k_{1,2}\in(-1,1)$ 
is equivalent to the condition 
$k^*_{1,2}>0$. This is a necessary and sufficient condition 
for the foliation by surfaces $S_\rho$
to extend throughout $M$. If this condition is satisfied the map $\Eps_\rho$ is 
a homeomorphism for any $\rho$.
\end{proof}

Interestingly, this condition makes sense not only in the quasi-Fuchsian situation
but is more general. Thus, for example, it applies to Schottky manifolds. But for
Schottky manifolds, with their single boundary component, the foliation by equidistant surfaces
$S_\rho$ cannot be smooth for arbitrary $\rho$. It is clear that surfaces must develop singularities
for some value of $\rho$. We therefore get an interesting corollary: 
\begin{cor}
There is no Liouville field $\phi$ on $\C$ invariant under a Schottky group such that $\phi_{z\bar{z}}$ is greater
than $|\phi_{zz}-(1/2)\phi_z^2|$ everywhere on $\C$. 
\end{cor}
\begin{proof}Indeed, if such a Liouville field existed,
we could have used it to construct a smooth equidistant foliation for arbitrary values of $\rho$,
but this is impossible. 
\end{proof}
A similar statement holds for a Kleinian group with more than two components
of the domain of discontinuity.

\section{The Schl\"afli formula ``from infinity''\index{Schl\"afli formula}\index{formula!Schl\"afli}}

In this section we obtain a formula for the variation of the
renormalized volume. 

\subsection{The Schl\"afli formula for hyperbolic polyhedra} \label{ssc:schlafli}

Recall first the classical Schl\"afli formula
(see e.g. \cite{milnor-schlafli}), which is a good a motivation for what follows.
Consider a hyperbolic polyhedron $P$. Under a first-order deformation of 
$P$, the first-order variation of the volume of $P$ is given by:
\begin{equation}\label{eq:schlafli}
dV = \frac{1}{2} \sum_e L_e d\theta_e~. 
\end{equation}
Here the sum is over the edges of $P$, $L_e$ is the length of the
edge $e$, and $\theta_e$ is its {\it exterior} dihedral angle. 

There is also an interesting ``dual'' Schl\"afli formula.
Let 
$$ V^* = V - \frac{1}{2}\sum_e L_e\theta_e~, $$
be the {\it dual volume} of $P$, then, still under a first-order
deformation of $P$,
\begin{equation}\label{eq:schlafli-dual}
dV^* = - \frac{1}{2} \sum_e \theta_e dL_e~. 
\end{equation}
This follows from the Schl\"afli formula (\ref{eq:schlafli})
by an elementary computation.

\subsection{The Schl\"afli formula for hyperbolic manifolds with boundary}

As we have seen in the previous sections, the renormalized volume of a convex
co-compact hyperbolic 
3-manifold $M$ can be expressed as the W-volume of any convex domain $N\subset
M$. The W-volume is equal to the 
volume of $N$ minus the quarter of the integral of the mean curvature over the
boundary of $N$.  
Let us consider what happens if one changes the metric in $M$.
As was shown in \cite{sem-era}, the following formula for the variation of the volume holds
\be
2 \delta V(N) = \int_{\partial N} 
\left( \delta H + \frac{1}{2} \langle\delta I,\II\rangle  \right) da.
\ee
Here $H$ is the trace of the shape operator $B=I^{-1} \II$, and the expression
$\langle A,B\rangle $ stands for
$\tr(I^{-1} A I^{-1} B)$. We can use this to get the
following expression for the variation of the W-volume:
\be\nonumber
\delta W(N) = \frac{1}{2} \int_{\partial N} \left( \delta H + \frac{1}{2} \langle \delta I,\II\rangle  \right) da - 
\frac{1}{4} \int_{\partial N} \delta H da - \frac{1}{4} \int_{\partial N} H
\delta(da)~, \ee
so that
\be \label{sch-1}
\delta W(N) = 
\frac{1}{4} \int_{\partial N} \left( \delta H + \langle \delta I,\II- \frac{H}{2} I\rangle  \right) da~.
\ee
To get the last equality we have used the obvious equality
\be
\delta da = \frac{1}{2} \tr(I^{-1} \delta I) = \frac{1}{2} \langle \delta I,I\rangle  da.
\ee
The formula (\ref{sch-1}) can be further modified using
\be
\delta H = \delta(\tr(I^{-1} \II)) = -\tr( I^{-1} (\delta I) I^{-1} \II) + \tr(I^{-1} \delta\II) = 
\ee
$$ = -\langle \delta I,\II\rangle  + \langle I,\delta \II\rangle~. $$
We get
\be\label{sch-2}
\delta W(N) = \frac{1}{4} \int_{\partial N} \langle \delta \II - \frac{H}{2} \delta I, I\rangle  da~.
\ee
It is this formula that will be our starting point for transformations to express the
variation in terms of the data at infinity.

\subsection{Parametrization by the data at infinity}

Let us now recall that given the data $I,\II$ on the boundary of $N$ one can
introduce the 
first and second fundamental forms ``at infinity'' via
(\ref{ff*}). Conversely, knowing 
the fundamental forms $I^*,\II^*$ ``at infinity'' one can recover the
fundamental forms 
on $\partial N$ via (\ref{metric-rel}). Our aim is to rewrite the variation
(\ref{sch-2}) 
of the W-volume in terms of the variations of the forms $I^*,\II^*$. 

\begin{lemma} \label{lm:schafli-inf}
The first-order variation of $W$ can be expressed as
\be\label{dw}
\delta W(N) = - \frac{1}{4} \int_{\partial N} \langle  \delta \II^* -
\frac{H^*}{2} \delta I^*, I^*\rangle  da^*. 
\ee
\end{lemma}

A proof can be found in \cite{volume}, it follows from a direct computation
based on the formulas expressing $I,\II$ in terms of $I^*,\II^*$. 

This formula could be compared to a more general formula applicable in any dimension given in
\cite{skenderis-solod} and reviewed in \cite{anderson:geometric}. A similar formula for the case of
4-dimensional manifolds is given by Anderson in \cite{anderson-L2}, and more recently by 
Albin \cite{albin} in higher dimensions. Our derivation here is different from that in the
cited references, for we interpret the variational formula as a version of the Schl\"afli formula.
We also note that the situation is simpler for even-dimensional manifolds, since
the renormalized volume is then canonically defined, while in odd
dimension it depends on the choice of a metric in the conformal class at
infinity. In odd dimension there is another, canonically defined
``renormalized volume'', namely the logarithmic term in the asymptotic expansion
of the volume as a function of the parameter of an equidistant foliation. This
quantity is different from the one used here (defined in (\ref{rv}) which is
the constant term in the same asymptotic expansion.

Formula (\ref{dw}) looks very much like the original formula (\ref{sch-2}),
except for the minus sign 
and the fact that the quantities at infinity are used. The fact that we have
got the same variational 
formula as in terms of the data on $\partial N$ is not too surprising. Indeed,
the variational 
formula (\ref{dw}) was obtained from (\ref{sch-2}) by applying the
transformation (\ref{metric-rel}). As it is clear from  
(\ref{ff*}), this transformation applied twice gives the identity map. In view
of this, it is hard to think of any other 
possibility for the variational formula in terms of $\delta I^*, \delta \II^*$
except being given by the same 
expression (\ref{sch-2}), apart from maybe with a different sign. This is
exactly what we see in (\ref{dw}).  

There is another expression of the first-order variation of $W$, dual to
(\ref{sch-1}), which will be useful below.

\begin{cor} \label{cr:variation}
The first-order variation of $W$ can also be expressed as
$$ \delta W = -\frac{1}{4}\int_{\dr N} \delta H^* + \langle \delta I^*,
\II^*_0\rangle da^*~, $$
where $\II^*_0$ is the traceless part (for $I^*$) of $\II^*$.
\end{cor}

\subsection{Conformal variations of the metric at infinity}

We can now use Corollary \ref{cr:variation}
to show that, when varying the W-volume with the area of the boundary
defined by the $I^*$ metric 
kept fixed, the variational principle forces the metric $I^*$ to have
constant negative curvature. 
The variations we consider here do not change the conformal
structure of the metric $I^*$, 
and thus do not change the manifold $M$. Geometrically they correspond to
small movements of the 
surface $\partial N$ inside the fixed manifold $M$. 
Thus, let us consider a conformal deformation of the metric $I^*$ of the form 
$\delta I^*=2u I^*$, where $u$ is some function on $\dr N$. Clearly for 
such variations $\langle \delta I^*,\II^*_0\rangle=0$, precisely because
$\II^*_0$ is traceless.

Let us consider the following functional
\be\label{var-3}
F(N) = W(N) - \frac{\lambda}{4} \int_{\partial N} da^* 
\ee
appropriate for finding an extremum of the W-volume with the area computed
using the 
metric $I^*$ kept fixed. The first variation of this functional gives,
using Corollary \ref{cr:variation}:
$$ \delta F = -\frac{1}{4}\int_{\dr N} (\delta H^*) da^* - \frac{\lambda}{4}
\int_{\dr N} 4u da^* = \frac{1}{4}\int_{\dr N} (\delta K^*) da^* -
\frac{\lambda}{4} \int_{\dr N} 4u da^*~. $$
But 
$$ \delta\int_{\dr N}K^* da^*= \int_{\dr N} (\delta K^*) + 4u K^* da^* =0 $$
by the Gauss-Bonnet formula, so that
$$ \delta F = \int_{\dr N} (- uK^* - u\lambda) da^*~. $$
It follows that critical points of $F$ are characterized by the fact that
$K^*=-\lambda$. 
It is not hard to compute the second variation and show that 
the critical points of $F$ are local maxima, see \cite{volume} for
details.

\begin{remark} {\rm As we have already discussed, the renormalized volume $W(N)$
is actually a functional of metrics $I^*$ on all the boundary components of $M$.
We have just established that this functional has an extremum, for
variations that keep the total area of the boundary components fixed, at
the constant curvature metric $I^*$. However, this is precisely
the defining property of the Liouville functional we have 
discussed in the Introduction. This establishes the renormalized volume ---
Liouville action functional relation. Moreover, one can turn the
argument around and use the renormalized volume $W(N)$ (as a functional
of the metrics $I^*$ on all the boundary components) as a 
{\it definition} of the Liouville functional. This point of view
explains why there is not one, but a whole set of Liouville functionals 
--- depending on which hyperbolic three-manifold is used --- and it
also explains why it is so hard to define the Liouville functional
in intrinsically two-dimensional terms --- because it is in fact
a three-dimensional quantity.}
\end{remark}

\subsection{The renormalized volume as a function on Teichm\"uller space}

Let us now consider the renormalized volume as a function over the
Teichm\"uller space of $\dr N$; in other terms, for each conformal class on
$\dr N$, we consider the extremum of $W$ over metrics of given area within this
conformal class. We have just seen that this extremum is
obtained at the (unique) constant curvature metric. The main goal here
is to recover by simple differential geometric methods important results of
Takhtajan and Zograf \cite{TZ-Schottky}, 
Takhtajan and Teo \cite{Takhtajan:2002cc} -- 
showing that the renormalized volume provides a K\"ahler
potential for the Weil-Petersson metric. So the ``volume'' that we consider
here is now defined as follows.

\begin{df}
Let $g$ be a convex co-compact hyperbolic metric on $M$, and let $c\in
\cT_{\dr M}$ be the conformal structure induced on $\dr_\infty M$. We call
$W_M(c)$ the value of $W$ on the equidistant foliation of $M$ near infinity
for which $I^*$ has constant curvature $-1$.
\end{df}

In other terms, by the results obtained above, $W_M(c)$ is
the maximum of $W$ over the metrics at infinity which have the same area as a hyperbolic metric, for
each boundary component of $M$. Throughout this section
the metric at infinity $I^*$ that we consider is the hyperbolic metric, while
the second fundamental form at infinity, $\II^*$, is uniquely determined by
the choice that $I^*$ is hyperbolic. Its traceless part is denoted by
$\II^*_0$. 

\subsection{The second fundamental form at infinity as the real part of a holomorphic
quadratic differential}

It is interesting to remark that, in the context considered here -- when $I^*$
has constant curvature -- the second fundamental form at infinity has a
complex interpretation. This can be compared with the same phenomenon,
discovered by Hopf \cite{hopf}, for the
second fundamental form of constant mean curvature surfaces in 3-dimensional
constant curvature spaces.

\begin{lemma}
When $K^*$ is constant, $\II^*_0$ is the real part of a holomorphic quadratic 
differential (for the complex structure associated to $I^*$) on $\dr_\infty
M$.  This holomorphic quadratic differential is given explicitly by
(\ref{theta}).
\end{lemma}

\begin{proof}
By definition $\II^*_0$ is traceless, which means that it is at each point the
real part of a quadratic differential: $\II^*_0=Re(h)$. Moreover, we have seen 
in Remark \ref{rk:codazzi*} that $B^*$ satisfies the Codazzi equation,
$d^{\nabla^*} B^*=0$. It follows as for constant mean curvature surfaces (see
e.g. \cite{minsurf}) that $h$ is holomorphic relative to the complex structure
of $I^*$.
\end{proof}

\subsection{The second fundamental form as a Schwarzian derivative}

The next step is to show that, for each boundary component $\dr_iM$ of $M$, 
$\II^*_{0i}$ is just the real part of the Schwarzian derivative of a natural equivariant 
map from the hyperbolic plane (with its canonical
complex projective structure) to $\dr_iM$ with its complex projective structure
induced by the hyperbolic metric on $M$. In the terminology used by 
McMullen \cite{McMullen}, $\II^*_{0i}$ is the difference between the complex
projective structure at infinity on $\dr_iM$ and the Fuchsian projective
structure on $\dr_iM$.

To state this result, let us 
call $\sigma_F$ the ``Fuchsian'' complex projective structure on $\dr_iM$,
obtained by applying the Poincar\'e uniformization theorem to the conformal
metric at infinity on $\dr_iM$. The universal cover of $\dr_iM$,
with the complex projective structure lifted from $\sigma_F$, is projectively
equivalent to a disk in $\C P^1$. We also call $\sigma_{QF}$ the projective
structure induced on $\dr_iM$ by the hyperbolic metric on $M$. Here ``$QF$''
stands for quasi-Fuchsian (while $M$ is only supposed to be convex
co-compact). This notation is used to keep close to the notation in
\cite{McMullen}. The map $\phi:(\dr_iM,\sigma_F) \rightarrow
(\dr_iM,\sigma_{QF})$ is conformal but not projective between $(\dr_iM,\sigma_F)$
and $(\dr_iM,\sigma_{QF})$, so we can consider its Schwarzian derivative 
$\cS(\phi)$. 

\begin{lemma} \label{lm:schwarzian}
$\II^*_0=-Re(\cS(\phi))$. 
\end{lemma}

A simple way to prove this assertion
is to use the formula (\ref{theta}) for the holomorphic quadratic differential
$\theta$ whose real part gives the traceless part of $\II^*$. The Liouville
field $\phi$ that enters into this formula can be simply expressed in
terms of the conformal map from $\dr_iM$ to the hyperbolic plane. It is
then a standard and simple computation to verify that $\theta$ is equal to
the Schwarzian derivative of this map, see e.g. \cite{TZ-Schottky}. 

It is possible to reformulate the statement (\ref{lm:schwarzian}) slightly 
using the following definition.

\begin{df}
Let $\theta_i$ be the restriction to $\dr_i M$ of $\cS(\phi)$.
\end{df}

This is notation is analogous to the notations used in \cite{McMullen}, the
index $i$ is useful to recall that this quantity is related to $\dr_iM$. 

Then $\theta_i$ is a quadratic holomorphic differential (QHD) on $\dr_iM$,
and, still using the notations in \cite{McMullen}, the definition of
$\theta_i$ can be rephrased as: $\theta_i=\sigma_{QF}-\sigma_F$. 
The Lemma can then be written as: $\II^*_{0i}=Re(\theta_i)$. 
A geometric proof of this lemma is given in the appendix of \cite{volume}.

\begin{remark} 
{\rm Note that
$\theta_i$ can also be considered as a complex-valued 1-form on the
Teichm\"uller space of $\dr_iM$. Indeed, it is well known that the cotangent
vectors to $\cT_S$, where $S$ is a Riemann surface, can be described as 
holomorphic quadratic differentials $q$ on $S$. The pairing with a tangent vector (Beltrami differential $\mu$) is given
by the integral of $q\mu$ over $S$. The complex structure on $\cT_S$ can then be described
as follows: the image of the cotangent vector $q$ under the action of the complex structure $J$ is
simply $J(q)=iq$. Another, more geometric way to state the action of $J$ is to note that
it exchanges the horizontal and vertical trajectories of $q$. 
Thus, holomorphic quadratic differentials $q$ on $S$ are actually
holomorphic 1-forms on $\cT_S$.}
\end{remark}

\subsection{The second fundamental form as the differential of $W_M$}

There is another simple interpretation of the traceless part of the
second fundamental form at infinity.

\begin{lemma} \label{lm:gradient}
The differential $dW_M$ of the renormalized volume $W_M$, 
as a 1-form over the Teichm\"uller space of $\dr M$, 
is equal to $(-1/4)\II^*_0$.
\end{lemma}

\begin{proof}
This is another direct consequence of Corollary \ref{cr:variation} 
because, as one
varies $I^*$ among hyperbolic metrics, $H^*$ (which is equal to $K^*$) remains
equal to $-1$, so that $\delta H^*=0$.  
\end{proof}

\begin{cor}
$\theta_i=-4\partial W_M$.
\end{cor}

\begin{proof}
This follows directly from the lemma, since we already know that $\theta_i$ is
a holomorphic differential.  
\end{proof}

\begin{remark} {\rm 
We would like to emphasize how much simpler is the proof given above
than that given in \cite{TZ-Schottky}, \cite{Takhtajan:2002cc}. Unlike in these
references, which obtain the above result on the gradient of $W_M$ via an
involved computation using a reasonably complicated cohomology
machinery, Corollary \ref{cr:variation} implies this result in
one line. This demonstrates the strength of the geometric method used here. Our
proof can be immediately extended even to situations where the methods of
\cite{Takhtajan:2002cc} are inapplicable, such as manifolds with cone
singularities. See more remarks on this case below.}
\end{remark}

Lemma \ref{lm:gradient} and, in particular, its corollary above,
is the key fact needed to demonstrate that the renormalized volume plays
the role of the K\"ahler potential on Teichm\"uller space. The remainder
of the proof is in part \ref{ssc:kahler} below, after some considerations on quasifuchsian
reciprocity, which are partly based on the results of this section and are needed
in the proof.

\section{Kleinian reciprocity\index{reciprocity!Kleinian}}

\subsection{Statement}

Kleinian reciprocity, as defined by McMullen \cite{McMullen}, is the extension of 
Theorem \ref{tm:mcm2} to the more general setting of a geometrically finite hyperbolic
3-manifold $M$. Let $\GF(M)$ be the space of complete geometrically finite 
hyperbolic metrics on $M$. Each metric $g\in \GF$ induces a complex projective
structure $\sigma(g)$ on $\dr M$. 

\begin{thm} \label{tm:mcm3}
$\sigma(\GF)$ is Lagrangian in $\CP_{\dr M}$.
\end{thm}

As explained in the introduction and proved in Section \ref{sssc:from} below, 
Theorem \ref{tm:mcm} is a direct corollary of this
statement.
The proof of Theorem \ref{tm:mcm}
given in \cite{McMullen} can be described as analytic, as it takes place in
the universal cover of $M$ and uses the symmetry of a certain kernel related to the
Beltrami problem. By contrast, the arguments considered here are mostly geometric and
take place in $M$.

We will describe here three (other) proofs of Theorem \ref{tm:mcm3}, corresponding
to different ways to consider the space of complex projective structures $\CP$.
\begin{itemize}
\item When $\CP$ is considered in complex terms, and identified with 
$T^*\cT_C$, the cotangent bundle of the space
of complex structures on $\dr M$, Theorem \ref{tm:mcm3} follows from the first
variation of the renormalized volume, as explained in the introduction. This is the
argument described in the introduction. 
\item When $\CP$ is considered in hyperbolic terms, and identified with $T^*\cT_H$,
Theorem \ref{tm:mcm3} follows from the dual Bonahon-Schl\"afli formula for the
first variation of the dual volume of the convex core. The equivalence with the
previous viewpoint is clear through Theorem \ref{tm:cp}.
\item When $\CP$ is considered as (a connected component of) a space of equivalence
classes of representations of $\pi_1(\dr M)$ into $PSL(2,\C)$, Theorem \ref{tm:mcm3}
can be proved by a completely different argument, based on exact sequences and Poincar\'e
duality, which was discovered (previousy) by S. Kerckhoff. The equivalence with the 
complex or the hyperbolic viewpoint follows from the fact that the Goldman form on
the space of representations of $\pi_1(\dr M)$ into $PSL(2,\C)$ is equal (up to
multiplication by a constant) to the symplectic form obtained from the cotangent
symplectic form on $T^*\cT$, as proved by Kawai \cite{kawai}.
\end{itemize}

We briefly describe these three arguments in the next subsections. We consider here for
simplicity the case of quasifuchsian manifolds, however all arguments can be
extended without difficulty to the more general context of geometrically finite
hyperbolic 3-manifolds.

\subsection{The first variation of the renormalized volume} \label{ssc:renormalized}

We give here the proofs announced in the introduction of quasifuchsian reciprocity
from the first-order variation formula for the renormalized volume. 

\subsubsection{From Theorem \ref{tm:mcm2} to Theorem \ref{tm:mcm}.} 
\label{sssc:from}

Suppose that Theorem \ref{tm:mcm2} holds, and consider a pair $(c'_-,0), (0,c'_+)$ 
of tangent vectors in $T_{c_-}\cTb\times T_{c_+}\cT$. In addition to the notations 
from Section \ref{ssc:19}, we introduce the notation $B:\cT_{\dr M}\rightarrow 
T^*\cT_{\dr M}$ for the map to the total space of the bundle $T^*\cT_{\dr M}$, 
corresponding to the section $\beta$ of $T^*\cT_{\dr M}$. We use the notation
$D$ for the differential as in the introduction, and call $\nabla$ the Levi-Civita
connection of the Weil-Petersson metric on the cotangent bundle of both $\cT$ and 
$\cTb$ (we could use another connection). Thus $\nabla$ determines a connection 
on $\cT_ {\dr M}=\cT\times \cTb$, and this connection defines an identification 
$T(T^*\cT_{\dr M})\simeq T\cT_{\dr M}\times T^*\cT_{\dr M}$. With those notations:
\begin{eqnarray*}
  \langle D\beta_+(c'_-,0), c'_+\rangle & - & \langle D\beta_-(0,c'_+), c'_-\rangle \\
& = &  \langle (\nabla_{c'_-}\beta_-, D\beta_+(c'_-,0)), (0,c'_+)\rangle \\
& - & \langle (D\beta(0,c'_+),\nabla_{c'_+}\beta_+), (c'_-,0)\rangle \\
& = &  \omega_*(B_*(c'_-,0), B_*(0,c'_+))= 0~,
\end{eqnarray*}
where the last equality follows from Theorem \ref{tm:mcm2}. 
This proves Theorem \ref{tm:mcm}.

Corollary \ref{cr:closed} is also a direct consequence of Theorem \ref{tm:mcm2}.
Indeed, for fixed $c_-$, $\beta_+(c_-,\cT)$ is a Lagrangian submanifold in $T^*\cT$ 
by Theorem \ref{tm:mcm2}. It is also the graph of a 1-form by definition, and
it is well-known (and can be checked by a direct computation) that a 1-form is
closed if and only if its graph is Lagrangian.

Conversely, Theorem \ref{tm:mcm2} follows quite directly from Theorem \ref{tm:mcm}
together with Corollary \ref{cr:closed}. To see this, let $c_-\in \overline{\cT},c_+\in
\cT, v_-\in T_{c_-}\overline{\cT}, v_+, v'_+\in T_{c_+}\cT$. Then the computation
above, done backwards, shows that
$$ \omega_*(B_*(v_-,0),B_*(0,v_+))=0~. $$
In addition, 
\begin{eqnarray*}
\omega_*(B_*(0,v_+),B_*(0,v'_+)) & = & 
\langle (D_{v_+}\beta_-,\nabla_{v_+}\beta_+), (0,v'_+)\rangle \\
& - & \langle (D_{v'_+}\beta_-,\nabla_{v'_+}\beta_+), (0,v_+)\rangle \\
& = & \langle \nabla_{v_+}\beta_+, v'_+\rangle
- \langle \nabla_{v'_+}\beta_+, v_+\rangle \\
& = & d\beta_+(v_+,v'_+) \\
& = & 0~, 
\end{eqnarray*}
where the last equality comes from Corollary \ref{cr:closed}. Theorem \ref{tm:mcm2}
follows by linearity.

\subsubsection{Proof of Theorem \ref{tm:mcm2}.}

Theorem \ref{tm:mcm2} follows very directly from the basic properties of the renormalized
volume $W_M$ as they are described above. Indeed we have seen that 
$$ \beta = -\II^*_0 = \frac{1}{4} dW_M~. $$
So $\beta$ is closed. 
This argument extends as it is to the more general setting of Theorem \ref{tm:mcm3}.

\subsection{The boundary of the convex core and the grafting map} \label{ssc:ccore}

The second argument leading to quasifuchsian reciprocity is also based on hyperbolic
geometry, and more specifically on the geometry of the convex core of quasifuchsian
3-manifolds. It rests on an extension of the classical Schl\"afli identity for
convex cores of quasifuchsian manifolds, found by Bonahon \cite{bonahon-variations},
which leads to an analog of Theorem \ref{tm:mcm3} where the renormalized volume
is replaced by the volume of the convex core, and the cotangent bundle of
Teichm\"uller space is considered in ``hyperbolic'' terms. 

\subsubsection{The convex core of quasifuchsian manifolds.}

Let $M$ be a quasifuchsian hyperbolic 3-manifold. $M$ contains a smallest non-empty 
geodesically convex subset, its {\it convex core} $C(M)$, which is compact and
homeomorphic to $M$. The boundary of $C(M)$ is therefore the disjoint union of
two copies of a surface $S$ of genus at least $2$, which we call the ``upper'' and
``lower'' boundary components of $C(M)$. Since $C(M)$ is a minimal convex subset, 
it has no extreme points, so $\dr C(M)$ is locally convex and ruled (there is a 
geodesic segment of the ambient manifold, contained in $C(M)$, containing each
point). It follows that the induced metric $m$ on $\dr C(M)$ is hyperbolic (it has
constant curvature $-1$), but $\dr C(M)$ is pleated along a measured 
geodesic lamination $l$. More details can be found in \cite{thurston-notes}.

\subsubsection{Measured laminations as cotangent vectors.}

When thinking of Teichm\"uller space in terms of hyperbolic metrics on surfaces,
it is natural to associate its cotangent bundle with measured laminations, rather
than with holomorphic quadratic differentials. The identification goes as follows. 
Let $l\in \ML$ be a measured lamination, and let $m\in \cT$ be a hyperbolic metric,
both on $S$. It is then possible (see \cite{bonahon-geodesic}) to define the length of $l$ for 
$m$, $L_m(l)$. This defines a smooth function 
$$ L_\cdot(l):\cT\rightarrow \R_{\geq 0}~. $$
The differential $dL_\cdot(l)$ at $m$ is an element of $T^*_m\cT$, and the
map $\ML\rightarrow T^*_m\cT$ is a homeomorphism (see e.g. \cite{cp}). 

This construction defines a natural identification between $\cT\times \ML$ 
and $T^*\cT$. But $T^*\cT$ has a cotangent symplectic structure, which we
call $\omega_H$ here. It can be pulled back to $\cT\times \ML$, where we still
call it $\omega_H$. It can be defined quite explicitly in terms of the 
intersection form on $\ML$, see \cite{sozen-bonahon}.

\subsubsection{A Lagrangian submanifold.}

Given a quasifuchsian metric $g\in \QF$, we can consider the induced metrics
on the upper and lower boundary components of the convex core, $m_+,m_-\in \cT$,
and the corresponding measured bending laminations, $l_+,l_-\in \ML$. So we 
have two points $(m_+,l_+),(m_-,l_-)\in \cT\times \ML$. This defines a map 
$H:\QF\rightarrow T^*\cT_{\dr M}$. 

\begin{thm} \label{tm:cc}
$H(\QF)$ is a Lagrangian submanifold of $(T^*\cT_{\dr M},\omega_H)$.
\end{thm}

The main ideas of the proof are explained in the next subsection. 
Theorem \ref{tm:mcm2} directly follows from this result and from Theorem
\ref{tm:cp}, according to which the grafting map is symplectic (up to 
a constant).

\subsubsection{The Bonahon-Schl\"afli formula for the volume of the convex core.}

The convex core of a quasifuchsian manifold is, in some ways, reminiscent
of a convex polyhedron. The main differences are: it has no vertices and
edges are replaced by a measured lamination. This gives, in a sense, a much
richer structure.

Bonahon \cite{bonahon} has extended the Schl\"afli formula recalled in
Part \ref{ssc:schlafli} to 
this setting as follows. Let $M$ be a convex co-compact hyperbolic
manifold (for instance, a quasifuchsian manifold), let $\mu$
be the induced metric on the boundary of the convex core, and let 
$\lambda$ be its measured bending lamination. By a ``first-order
variation'' of $M$ we mean a first-order variation of the 
representation of the fundamental group of $M$. Bonahon shows that
the first-order variation of $\lambda$  under a first-order
variation of $M$ is described by a 
transverse H\"older distribution $\lambda'$, and there is a 
well-defined notion of length of such transverse H\"older 
distributions. This leads to 
a version of the Schl\"afli formula.

\begin{lemma}[The Bonahon-Schl\"afli formula \cite{bonahon}]\label{lm:bs}
The first-order variation of the volume $V_C$ of the convex
core of $M$, under a first-order variation of $M$, is given by
$$ dV_C = \frac{1}{2} L_\mu(\lambda')~. $$
\end{lemma}

Here $\lambda'$ is the first-order variation of the measured bending lamination,
which is a H\"older cocycle so that its length for $\mu$ can be defined, see 
\cite{bonahon-toulouse,bonahon-ens,bonahon,bonahon-variations}. 

\subsubsection{The dual volume.}

Just as for polyhedra above, we define the dual volume
of the convex core of $M$ as 
$$ V_C^* = V_C -\frac{1}{2} L_\mu(\lambda)~. $$

\begin{lemma}[The dual Bonahon-Schl\"afli formula] \label{lm:bs-dual}
The first-order variation of $V^*$ under a first-order 
variation of $M$ is given by
$$ dV_C^* = - \frac{1}{2} L'_\mu(\lambda)~. $$
\end{lemma}

This formula has a very simple interpretation in terms of 
the geometry of Teichm\"uller space: up to the factor $-1/2$,
$dV^*$ is equal to the pull-back by $\delta$ of the Liouville form of the cotangent
bundle $T^*\cT$. Note also that this formula can be 
understood in an elementary way, without reference to a
transverse H\"older distribution: the measured lamination $\lambda$ is
fixed, and only the hyperbolic metric $\mu$ varies. 
The proof can be found in \cite{cp}, it is based on Lemma \ref{lm:bs}.

Theorem \ref{tm:cc} is a direct consequence of Lemma 
 \ref{lm:bs-dual}: since $dV_C^*$ coincides with the
Liouville form of $T^*\cT_{\dr M}$ on $H(\QF)$, it follows
immediately that $H(\QF)$ is Lagrangian for the symplectic
form $\omega_H$ on $T^*\cT_{\dr M}$.

\subsection{Deformations of representations and Poincar\'e duality} \label{ssc:kerckhoff}

Steve Kerckhoff found another (unpublished) proof of Theorem \ref{tm:mcm3} based
on topological ideas and in particular on his earlier work with Craig Hodgson \cite{hk}.

This proof works in the context of deformations of $PSL(2,\C)$ representations,
the symplectic form on $\CP$ is here the Goldman symplectic form $\omega_G$ on $\CP$.
Recall (from \cite{goldman-symplectic}) that given a complex projective structure 
$\sigma$ on $\dr M$, its holonomy representation is a morphism $\rho$ from $\pi_1(\dr M)$
to $PSL(2,\C)$. The tangent space to $\CP$ at $\sigma$ is then naturally identified
with the cohomology space $H^1(\dr M; E)$, where $E$ is a $sl(2,\C)$ bundle over
$\dr M$ naturally associated to $\rho$ -- it is the bundle of local projective 
vector fields for $\sigma$ on $\dr M$. Given two cohomology classes 
$u,v\in H^1(\dr M; E)$, one can consider their cup product $u\cup v\in H^2(\dr M,\C)$,
and integrate it over $\dr M$. This defines the Goldman symplectic form 
$$ \omega_G:H^1(\dr M; E)\times H^1(\dr M; E)\rightarrow \C~. $$

Kawai \cite{kawai} proved that this symplectic form is equal, up to a constant, to the
canonical symplectic form on $\CP_{\dr M}$, obtained for instance by identification 
of $\CP_{\dr M}$ with $T^*\cT_{\dr M}$ through the Schwarzian derivative, see 
\cite{dumas-survey}.

Kerckhoff's proof is based on the long exact sequence in cohomology for the pair
$(M,\dr M)$:
$$ \cdots \rightarrow H^1(M,\dr M; E) \rightarrow H^1(M; E)\stackrel{\alpha}{\rightarrow} 
H^1(\dr M; E)\stackrel{\beta}{\rightarrow} H^2(M,\dr M; E)\rightarrow \cdots $$
Here $\alpha$ is restriction of the deformation from $M$ to $\dr M$.
Note that the map $H^1(M,\dr M; E) \rightarrow H^1(M; E)$ is zero, since any non-trivial
deformation of the hyperbolic structure on $M$ induces a non-zero deformation of
the complex projective structure on the boundary (this follows for instance
from the Ahlfors-Bers theorem). As a consequence, $\alpha$ is injective.

Part of the long exact sequence above can be extended as the commutative
diagram below, taken from \cite[p. 42]{hk}.
$$
\xymatrix{
H^1(M; E)\ar[r]^\alpha\ar[d] & H^1(\dr M; E)\ar[r]^{\beta}\ar[d] & H^2(M,\dr M; E)\ar[d]\\
H^2(M,\dr M; E)^* \ar[r]^{\beta^*}& H^1(\dr M; E)^*\ar[r]^{\alpha^*} & H^1(M; E)^*
}
$$
Here the vertical arrows are the Poincar\'e duality maps. Recall that Poincar\'e
duality can be defined using the cup product as above. In particular, the Poincar\'e
dual $u^*$ of a form $u$, for instance in $H^1(\dr M; E)$, is such that, for all $v\in 
H^1(\dr M; E)$, $\omega_G(u,v)=\langle u^*, v\rangle$. 

Let $u,v\in H^1(M; E)$. It follows from the above commutative diagram that
\begin{eqnarray*}
\omega_G(\alpha(u),\alpha(v)) & = & \langle \alpha(u),\alpha(v)^*\rangle \\
& = & \langle \alpha(u), \beta^*(v^*)\rangle \\
& = & \langle \beta\circ\alpha (u), v^* \rangle \\
& = & 0~. 
\end{eqnarray*}
It also follows from the above exact sequence (or from the upper part of the
commutative diagram and the fact that $\alpha$ is injective) that 
$\dim H^1(\dr M; E)=2\dim H^1(M; E)$ (see \cite{hk} for the details). 
So $\alpha(H^1(M; E))$ is Lagrangian in $H^1(\dr M; E)$, and this, along with
Kawai's result \cite{kawai}, provides 
yet another proof of Theorem \ref{tm:mcm3}.

\section{The renormalized volume as a K\"ahler potential}
\label{sc:kaehler}

In this section we again consider the setting of quasi-Fuchsian manifolds and 
recover the result of Takhtajan and Teo \cite{Takhtajan:2002cc}: the renormalized volume $W_M$ with
$c_-$ fixed is a K\"ahler potential for the Weil-Petersson metric on
$\cT_{\dr_+M}$.

\subsection{Notations}

To simplify notations a little, we set
$\theta_{c_-}:=\theta(c_-,\cdot)$, so that the real part of 
$\theta_{c_-}$ is $\beta_+(c_-, \cdot)$. 
Since we already know that $\theta_{c_-}=4\partial W_M$, we
only need to prove that $\db(i\theta_{c_-}) = -2\omega_{WP}$, where
$\omega_{WP}$ is the K\"ahler form of the Weil-Petersson metric on
$\cT_{\dr_+M}$. 

An important part of the argument is that $d\theta_{c_-}$, as a 2-form on
$\cT_{\dr_+M}$, does not depend on $c_-$. This appears as Theorem 7.2 
in McMullen's paper \cite{McMullen}. We include a proof for completeness,
following the proof given in \cite{McMullen}.

\begin{prop}
The differential $d\theta_{c_-}$, considered as a complex-valued 2-form
on $\cT_{\dr_+M}$, does not depend on $c_-$.
\end{prop}

\begin{proof}
Let $v_-\in T_{c_-}\cT_{\dr_-M}$, we want to show that 
the corresponding first-order variation $D_{v_-}(d\theta_{c_-})$
of $d\theta_{c_-}$ vanishes. This will follow from the fact that 
the first-order variation of $\theta_{c_-}$ corresponding to $v_-$,
$D_{v_-}\theta_{c_-}$, is the differential of a function defined 
on $\cT_{\dr_+M}$, namely the function $f_{v_-}$ defined by
$$ f_{v_-}(c_+) = \langle \theta_{c_-}(c_+), v_-\rangle~, $$
where $\langle , \rangle$ is the duality pairing.

The fact that $D_{v_-}\theta_{c_-}=df_{v_-}$ can be proved by
evaluating both sides on a vector $v_+\in T_{c_+}\cT_{\dr_+M}$ and using
the quasi-Fuchsian reciprocity. Since $\theta_{c_-}$ is a 
complex 1-form with real part equal to $\beta_+(c_-,\cdot)$, Theorem
\ref{tm:mcm} indicates that:
$$ \langle D_{v_-}\theta_{c_-}, v_+\rangle 
= \langle D\beta_+(c_-, c_+)(v_-,0), v_+,\rangle = $$
$$ = \langle D\beta_-(c_-,c_+)(0,v_+), v_-\rangle = 
df_{v_-}(v_+)~. $$
It clearly follows that $d\theta_{c_-}$, as a 2-form on $\cT_{\dr_+M}$,
does not depend on $c_-$.
\end{proof}

\subsection{Local deformations near the Fuchsian locus}

That $W_M$ is a K\"ahler potential 
is now reduced to a simple computation in the Fuchsian
situation. 

\begin{prop} \label{pr:varII}
Suppose that $M$ is a Fuchsian manifold, with $c_+=c_-$. Let $I^*$ be the
hyperbolic metric in the conformal class $c_+$. Under a first-order
deformation which does not change $c_-$, 
the variation of $I^*$ and $\II^*_0$ on $\dr_+M$ are related by:
$$ \delta \II^*_0 = - \delta I^*~. $$
\end{prop}

The proof is quite elementary, it is based on the fact that, for a quasifuchsian
manifold which is ``close to Fuchsian'', the metrics at infinity on the upper
and lower components of the boundary are $I^*$ and $\III^*$ respectively (where
$\III^*$ is the third fundamental form at infinity on the upper boundary component). 
Moreover, if $I^*$ has constant curvature, then $\III^*$ also has constant
curvature (see \cite{volume}). 

\subsection{K\"ahler potential} \label{ssc:kahler}

We can reformulate this statement by calling $\theta_R:=Re(\theta_{c_-})$, so
that, by Lemma \ref{lm:schwarzian}, 
$\theta_R(X)=- \langle X,\II^*_0\rangle$. Using the previous proposition,
this can then be stated as
$$ (D_X\theta_R)(Y) = \langle X,Y\rangle_{WP}~, $$
where $D$ is the Levi-Civita connection of the Weil-Petersson metric on
$\cT_{\dr_+M}$. 

We can now compute explicitly an expression of $\db\theta_{c_-}$, denoting by
$J$ the complex structure on $\cT_{\dr_-M}$.
\begin{eqnarray*} 
\db\theta_{c_-}(X,Y) & = & (D_X\theta_{c_-})(Y) + i(D_{JX}\theta_{c_-})(Y) \\
& = & (D_X\theta_R)(Y) - i(D_X\theta_R)(JY) + \\
& & + i((D_{JX}\theta_R)(Y) -
i(D_{JX}\theta_R)(JY)) \\
& = & \langle X,Y\rangle - i \langle X,JY\rangle + i\langle JX,Y\rangle +
\langle JX,JY\rangle \\
& = & 2(\langle X,Y\rangle -i\langle X,JY\rangle)~. 
\end{eqnarray*}

This means precisely that $\db \theta_{c_-}(X,JX)=2i\|X\|^2_{WP}$, and we
recover the result of Takhtajan and Teo \cite{Takhtajan:2002cc} that $W_M$ is
a K\"ahler potential for the Weil-Petersson metric.

Note that this statement could be compared to the fact (proved recently in 
\cite{Guo-Huang-Wang}) that, on the Fuchsian locus, the Hessian of the 
area of the (unique) closed minimal surface, as a function on $\overline{\cT}
\times \cT$, is compatible with the hypothesis that this area is also a 
K\"alher potential for the Weil-Petersson metric.

\section{The relative volume of hyperbolic ends}
\label{sc:relative}

So far we have considered a version of the renormalized volume defined for
hyperbolic 3-manifolds $M$ as a whole. This means that only certain very special
projective structures can arise at boundary components $\partial M$. It is
interesting, however, to extend the notion of renormalized volume
(and thus of Liouville action) to arbitrary projective structures at
the boundary. This is achieved by the notion of the relative volume that
we consider in this section. When the projective structure in question is
such that a non-singular hyperbolic 3-manifold $M$ realizing it exists,
then the renormalized volume of $M$ is just the sum of relative volumes
of its hyperbolic ends and the (dual) volume of the convex core $C(M)$.

We first recall some results due to Bonahon 
concerning the first variation of the volume of the 
convex core of a quasifuchsian manifold. We then introduce
the relative volume of a hyperbolic end, and give a 
first variation formula for it. This establishes an analog
of Kleinian reciprocity in the relative volume context,
and proves that the grafting map is symplectic.

\subsection{Definition}
\label{ssec:relative}

The relative volume is defined for
(geometrically finite) hyperbolic ends rather than for hyperbolic manifolds. 
Thus, consider a hyperbolic end $M$. The procedure used in the 
definition of the renormalized volume can be used in this setting, 
leading to the relative volume of the end. We will say that a geodesically
convex subset $K\subset M$ is a {\it collar} if it is relatively compact and 
contains the metric boundary $\dr_0M$ of $M$ (possibly all geodesically 
convex relatively compact subsets of $M$ are collars, but it is not necessary
to consider this question here). Then $\dr K\cap M$ is a locally convex surface in $M$.

The relative volume of $M$ is related both to the (dual) volume of the convex
core and to the renormalized volume; it is defined as the renormalized volume,
but starting from the metric boundary of the hyperbolic end. We follow
the same path as for the renormalized volume and start from a collar $K\subset M$. We set 
$$ W(K) = V(K) - \frac{1}{4}\int_{\dr K} H da + 
\frac{1}{2} L_\mu(\lambda)~, $$
where $H$ is the mean curvature of the boundary of $K$, $\mu$ is the 
induced metric on the metric boundary $\partial_0 M$ of $M$, and $\lambda$ is its
measured bending lamination.

As for the renormalized volume we define the metric at infinity as
$$ I^* := \lim_{\rho\rightarrow \infty} 2e^{-2\rho}I_\rho~, $$
where $I_\rho$ is the set of points at distance $\rho$ from $K$.
The conformal structure of $I^*$ is equal to the canonical conformal 
structure $c_\infty$ at infinity of $M$.

Here again, $W$ only depends on $I^*$ (and on $M$). Not all metrics in $c_\infty$ can 
be obtained from a compact subset of $M$, however all metrics do define
an equidistant foliation close to infinity in $M$, and it is still possible
to define $W(I^*)$ even when $I^*$ is not obtained from a convex subset of
$M$. So $W$ defines a function,
still called $W$, from metrics in the conformal class $c_\infty$ to $\R$. 

\begin{lemma}
For fixed area of $I^*$, $W$ is maximal exactly when $I^*$ has
constant curvature.  
\end{lemma}

The proof follows directly from the arguments used in \cite{volume} (and reviewed in Section 7)
so we do not repeat it here. It takes place entirely on the
boundary at infinity so it makes no difference whether one considers a hyperbolic end or a geometrically
finite hyperbolic manifold.

\begin{df}
The relative volume $V_R$ of $M$ is  $W(I^*)$ when $I^*$ is the
hyperbolic metric in the conformal class at infinity on $M$.
\end{df}

\subsection{The first variation of the relative volume}

\begin{prop} \label{pr:relative}
Under a first-order variation of the hyperbolic end, the first-order variation
of the relative volume is given by
\begin{equation}
  \label{eq:s-rel}
  V'_R = \frac{1}{2} L'_\mu(\lambda) - \frac{1}{4} \int_{\dr_\infty E} 
\langle I^*{}',\II^*_0 \rangle da^*~.
\end{equation}
\end{prop}

The proof is based on the arguments described in the previous sections, both for the
first variation of the renormalized volume and for the first
variation of the volume of the convex core.

\subsection{The grafting map is symplectic}
\label{ssc:grafting}

Since hyperbolic ends are in one-to-one correspondence with $\C P^1$-structures,
we can consider the relative volume $V_R$ as a function on the space of projective structures $\CP$.
This space is canonically identified with the (complex) cotangent bundle 
$T^*\cT_C$, where the subscript $C$ indicates that one is talking about the 
complex Teichm\"uller space. Let $\beta_C$ be the Liouville form on $T^*\cT_C$.
Consider now the space $\cT\times \ML$ associated with the metric boundary of 
our hyperbolic end. In the previous sections we have discussed how this space
can be naturally identified by the map $\delta$ with $T^*\cT_H$. 
Let $\lambda_H$ be the Liouville form on $T^*\cT_H$.
 
We can now consider the grafting map $Gr: \cT\times\ML \to \C P^1$, and 
the composition $\delta\circ Gr^{-1}:\CP\rightarrow T^*\cT_H$. This latter map
turns out to be $C^1$ (see \cite{cp}), a fact which is somewhat surprising since
there is no $C^1$ structure on $\ML$. 
It pulls back $\lambda_H$ as
$$ (\delta\circ Gr^{-1})^* \lambda_H = L_\mu'(\lambda)~. $$
Under the identification of $\CP$ with $T^*\cT_C$ through the Schwarzian derivative,
the expression of $\lambda_C$ is
$$ \lambda_C = \int_{\dr_\infty M}\langle {I^*}',\II^*_0\rangle da^*~. $$
So Proposition \ref{pr:relative} can be formulated as
$$ dV_R = \frac{1}{2} (\delta\circ Gr^{-1})^*\lambda_H - \frac{1}{4} \lambda_C~, $$
and it follows that $2(\delta\circ Gr^{-1})^*\omega_H = \omega_C$. This means
that the grafting map preserves (up to a constant) the symplectic form and
is thus symplectic. This statement can also be rephrased in a way analogous to (\ref{tm:mcm3})
by saying that the subspace of the space $(\cT\times\ML) \times \CP$ that can be realized on the
two boundaries of a hyperbolic end is a Lagrangian submanifold in $(\cT\times\ML) \times \CP$.

\section{Manifolds with particles and the Teichm\"uller theory of surfaces
  with cone singularities} 
\label{sc:cone}

One key feature of the arguments presented in this work is that they are
always local, 
in the sense that they depend on local quantities defined on the boundaries of
compact subsets of quasi-Fuchsian manifolds. Thus, we make only a very limited
use 
of the fact that the quasi-Fuchsian manifolds are actually quotients of 
hyperbolic 3-space by a group of isometries. One place where this is used is
in the proof of the fact that $\II^*$ is determined by $I^*$ (actually a
direct consequence of the Bers double uniformization theorem).
We expect that all the results should extend from quasi-Fuchsian
(more generally geometrically finite) 
manifolds to the ``quasi-Fuchsian manifolds with particles''  
which were studied e.g. in \cite{minsurf,qfmp}. Those are actually
cone-manifolds, with 
cone singularities along infinite lines running from one connected component of
the boundary at infinity to the other, along which the cone angle is less than
$\pi$. 

In the (non-singular)
quasi-Fuchsian setting the Bers double uniformization theorem shows that
everything is determined by the conformal structure at infinity. 
The corresponding statement holds for ``quasifuchsian manifolds
with particles''; a first step towards it is made in \cite{qfmp}, while the
second step is made in \cite{conebend}. The corresponding statement for
more general, convex co-compact manifolds, remains however elusive.

Those results can be used to obtain
results on the Teichm\"uller-type space of hyperbolic metrics with $n$ cone
singularities of prescribed angles on a closed surface  of genus $g$. Note
that this space, which can be denoted by $\cT_{g,n,\theta}$ (with $\theta=
(\theta_1, \cdots, \theta_n)\in (0,\pi)^n$) is topologically the same as the
``usual'' Teichm\"uller space $\cT_{g,n}$ of hyperbolic metrics with $n$ cusps
(with a one-to-one correspondence from \cite{troyanov}) but it has a natural
``Weil-Petersson'' metric which is different. It follows from the
considerations made here, extended to quasifuchsian manifolds with particles,
that this ``Weil-Petersson'' metric is still K\"ahler, with the renormalized
volume playing the role of a K\"ahler potential --- a result also obtained
by different arguments by Schumacher and Trapani \cite{schumacher-trapani}. 
We leave the detailed investigation of this extension to quasifuchsian cone 
manifolds for future work.

\bibliographystyle{amsplain}

\def\cprime{$'$}
\providecommand{\bysame}{\leavevmode\hbox to3em{\hrulefill}\thinspace}
\providecommand{\MR}{\relax\ifhmode\unskip\space\fi MR }
\providecommand{\MRhref}[2]{%
  \href{http://www.ams.org/mathscinet-getitem?mr=#1}{#2}
}
\providecommand{\href}[2]{#2}

\end{document}